\documentclass[12pt]{article}
\usepackage{latexsym, amssymb}

\DeclareMathAlphabet\gothic{U}{euf}{m}{n}

\makeatletter
\def\eqnarray{\stepcounter{equation}\let\@currentlabel=\theequation
\global\@eqnswtrue
\tabskip\@centering\let\\=\@eqncr
$$\halign to \displaywidth\bgroup\hfil\global\@eqcnt\z@
  $\displaystyle\tabskip\z@{##}$&\global\@eqcnt\@ne
  \hfil$\displaystyle{{}##{}}$\hfil
  &\global\@eqcnt\tw@ $\displaystyle{##}$\hfil
  \tabskip\@centering&\llap{##}\tabskip\z@\cr}

\def\endeqnarray{\@@eqncr\egroup
      \global\advance\c@equation\m@ne$$\global\@ignoretrue}

\def\@yeqncr{\@ifnextchar [{\@xeqncr}{\@xeqncr[5pt]}}
\makeatother

\begin{document}
\bibliographystyle{tom}

\newtheorem{lemma}{Lemma}[section]
\newtheorem{thm}[lemma]{Theorem}
\newtheorem{cor}[lemma]{Corollary}
\newtheorem{voorb}[lemma]{Example}
\newtheorem{rem}[lemma]{Remark}
\newtheorem{prop}[lemma]{Proposition}
\newtheorem{stat}[lemma]{{\hspace{-5pt}}}
\newtheorem{obs}[lemma]{Observation}
\newtheorem{defin}[lemma]{Definition}

\newenvironment{remarkn}{\begin{rem} \rm}{\end{rem}}
\newenvironment{exam}{\begin{voorb} \rm}{\end{voorb}}
\newenvironment{defn}{\begin{defin} \rm}{\end{defin}}
\newenvironment{obsn}{\begin{obs} \rm}{\end{obs}}

\newenvironment{emphit}{\begin{itemize} }{\end{itemize}}

\newcommand{\gota}{\gothic{a}}
\newcommand{\gotb}{\gothic{b}}
\newcommand{\gotc}{\gothic{c}}
\newcommand{\gote}{\gothic{e}}
\newcommand{\gotf}{\gothic{f}}
\newcommand{\gotg}{\gothic{g}}
\newcommand{\gothh}{\gothic{h}}
\newcommand{\gotk}{\gothic{k}}
\newcommand{\gotm}{\gothic{m}}
\newcommand{\gotn}{\gothic{n}}
\newcommand{\gotp}{\gothic{p}}
\newcommand{\gotq}{\gothic{q}}
\newcommand{\gotr}{\gothic{r}}
\newcommand{\gots}{\gothic{s}}
\newcommand{\gotu}{\gothic{u}}
\newcommand{\gotv}{\gothic{v}}
\newcommand{\gotw}{\gothic{w}}
\newcommand{\gotz}{\gothic{z}}
\newcommand{\gotA}{\gothic{A}}
\newcommand{\gotB}{\gothic{B}}
\newcommand{\gotG}{\gothic{G}}
\newcommand{\gotL}{\gothic{L}}
\newcommand{\gotS}{\gothic{S}}
\newcommand{\gotT}{\gothic{T}}

\newcommand{\mn}{\marginpar{\hspace{1cm}*} }
\newcommand{\mnn}{\marginpar{\hspace{1cm}**} }

\newcommand{\mnq}{\marginpar{\hspace{1cm}*???} }
\newcommand{\mnnq}{\marginpar{\hspace{1cm}**???} }

\newcounter{teller}
\renewcommand{\theteller}{\Roman{teller}}
\newenvironment{tabel}{\begin{list}%
{\rm \bf \Roman{teller}.\hfill}{\usecounter{teller} \leftmargin=1.1cm
\labelwidth=1.1cm \labelsep=0cm \parsep=0cm}
                      }{\end{list}}

\newcounter{tellerr}
\renewcommand{\thetellerr}{(\roman{tellerr})}
\newenvironment{subtabel}{\begin{list}%
{\rm  (\roman{tellerr})\hfill}{\usecounter{tellerr} \leftmargin=1.1cm
\labelwidth=1.1cm \labelsep=0cm \parsep=0cm}
                         }{\end{list}}
\newenvironment{ssubtabel}{\begin{list}%
{\rm  (\roman{tellerr})\hfill}{\usecounter{tellerr} \leftmargin=1.1cm
\labelwidth=1.1cm \labelsep=0cm \parsep=0cm \topsep=1.5mm}
                         }{\end{list}}

\newcommand{\Ni}{{\bf N}}
\newcommand{\Ri}{{\bf R}}
\newcommand{\Ci}{{\bf C}}
\newcommand{\Si}{{\bf S}}
\newcommand{\Ti}{{\bf T}}
\newcommand{\Zi}{{\bf Z}}
\newcommand{\Fi}{{\bf F}}

\newcommand{\proof}{\mbox{\bf Proof} \hspace{5pt}} 
\newcommand{\remark}{\mbox{\bf Remark} \hspace{5pt}}
\newcommand{\ruimte}{\vskip10.0pt plus 4.0pt minus 6.0pt}

\newcommand{\simh}{{\stackrel{{\rm cap}}{\sim}}}
\newcommand{\ad}{{\mathop{\rm ad}}}
\newcommand{\Ad}{{\mathop{\rm Ad}}}
\newcommand{\Aut}{\mathop{\rm Aut}}
\newcommand{\arccot}{\mathop{\rm arccot}}
\newcommand{\capp}{{\mathop{\rm cap}}}
\newcommand{\rcapp}{{\mathop{\rm rcap}}}
\newcommand{\Capp}{{\mathop{\rm Cap}}}
\newcommand{\diam}{\mathop{\rm diam}}
\newcommand{\divv}{\mathop{\rm div}}
\newcommand{\dist}{\mathop{\rm dist}}
\newcommand{\codim}{\mathop{\rm codim}}
\newcommand{\RRe}{\mathop{\rm Re}}
\newcommand{\IIm}{\mathop{\rm Im}}
\newcommand{\Tr}{{\mathop{\rm Tr}}}
\newcommand{\Vol}{{\mathop{\rm Vol}}}
\newcommand{\card}{{\mathop{\rm card}}}
\newcommand{\supp}{\mathop{\rm supp}}
\newcommand{\sgn}{\mathop{\rm sgn}}
\newcommand{\essinf}{\mathop{\rm ess\,inf}}
\newcommand{\esssup}{\mathop{\rm ess\,sup}}
\newcommand{\Int}{\mathop{\rm Int}}
\newcommand{\Leibniz}{\mathop{\rm Leibniz}}
\newcommand{\lcm}{\mathop{\rm lcm}}
\newcommand{\loc}{{\rm loc}}
\newcommand{\Jac}{{\rm Jac}}

\newcommand{\mod}{\mathop{\rm mod}}
\newcommand{\spann}{\mathop{\rm span}}
\newcommand{\one}{1\hspace{-4.5pt}1}

\newcommand{\DWR}{}

\hyphenation{groups}
\hyphenation{unitary}

\newcommand{\tfrac}[2]{{\textstyle \frac{#1}{#2}}}

\newcommand{\cb}{{\cal B}}
\newcommand{\cc}{{\cal C}}
\newcommand{\cd}{{\cal D}}
\newcommand{\ce}{{\cal E}}
\newcommand{\cf}{{\cal F}}
\newcommand{\ch}{{\cal H}}
\newcommand{\ci}{{\cal I}}
\newcommand{\ck}{{\cal K}}
\newcommand{\cl}{{\cal L}}
\newcommand{\cm}{{\cal M}}
\newcommand{\cn}{{\cal N}}
\newcommand{\cp}{{\cal P}}
\newcommand{\co}{{\cal O}}
\newcommand{\cs}{{\cal S}}
\newcommand{\ct}{{\cal T}}
\newcommand{\cu}{{\cal U}}
\newcommand{\cx}{{\cal X}}
\newcommand{\cy}{{\cal Y}}
\newcommand{\cz}{{\cal Z}}

\newcommand{\Tan}{{\rm Tan}}

\newcommand{\wtozp}{W^{1,2}\raisebox{10pt}[0pt][0pt]{\makebox[0pt]{\hspace{-34pt}$\scriptstyle\circ$}}}
\newlength{\hightcharacter}
\newlength{\widthcharacter}
\newcommand{\covsup}[1]{\settowidth{\widthcharacter}{$#1$}\addtolength{\widthcharacter}{-0.15em}\settoheight{\hightcharacter}{$#1$}\addtolength{\hightcharacter}{0.1ex}#1\raisebox{\hightcharacter}[0pt][0pt]{\makebox[0pt]{\hspace{-\widthcharacter}$\scriptstyle\circ$}}}
\newcommand{\cov}[1]{\settowidth{\widthcharacter}{$#1$}\addtolength{\widthcharacter}{-0.15em}\settoheight{\hightcharacter}{$#1$}\addtolength{\hightcharacter}{0.1ex}#1\raisebox{\hightcharacter}{\makebox[0pt]{\hspace{-\widthcharacter}$\scriptstyle\circ$}}}
\newcommand{\scov}[1]{\settowidth{\widthcharacter}{$#1$}\addtolength{\widthcharacter}{-0.15em}\settoheight{\hightcharacter}{$#1$}\addtolength{\hightcharacter}{0.1ex}#1\raisebox{0.7\hightcharacter}{\makebox[0pt]{\hspace{-\widthcharacter}$\scriptstyle\circ$}}}

\newpage

 \begin{center}
  {\Large{\bf The weighted  Hardy constant  }} \\[5mm]
\large Derek W. Robinson$^\dag$ \\[1mm]

\normalsize{30th March 2021}\\[1mm]
\end{center}

\vspace{+5mm}

\begin{list}{}{\leftmargin=1.7cm \rightmargin=1.7cm \listparindent=15mm 
   \parsep=0pt}
   \item
{\bf Abstract} $\;$ 
Let $\Omega$ be a domain in $\Ri^d$ and $d_\Gamma$  the Euclidean distance to the boundary $\Gamma$.
We investigate whether the weighted Hardy inequality
\[
\|d_\Gamma^{\,\delta/2-1}\varphi\|_2\leq a_\delta\,\|d_\Gamma^{\,\delta/2}\,(\nabla\!\varphi)\|_2
\]
is valid, with $\delta\geq 0$ and $a_\delta>0$, for all $\varphi\in C_c^1(\Gamma_{\!\!r})$ and all small $r>0$ where $\Gamma_{\!\!r}=\{x\in\Omega: d_\Gamma(x)<r\}$.
First we prove that if $\delta\in[0,2\rangle$ then  the inequality  is equivalent to the weighted version of Davies' weak Hardy inequality  on $\Omega$
with equality of the corresponding optimal constants.
Secondly, we establish that if $\Omega$ is a uniform domain with Ahlfors regular boundary then 
 the inequality is satisfied for all $\delta\geq 0$, and all small $r$,  with the exception of the value $\delta=2-(d-d_{\!H})$ where $d_{\!H}$ is the Hausdorff dimension of $\Gamma$.
Moreover, the optimal constant $a_\delta(\Gamma)$ satisfies $a_\delta(\Gamma)\geq 2/|(d-d_{\!H})+\delta-2|$.
Thirdly,  if $\Omega$ is a $C^{1,1}$-domain or  a convex domain $a_\delta(\Gamma)=2/|\,\delta-1|$ for all $\delta\geq0$ with $\delta\neq1$.
The same conclusion is correct if $\Omega$ is the complement of a convex domain and $\delta>1$ but if  $\delta\in[0,1\rangle$ then  $a_\delta(\Gamma)$ can be strictly larger  than $2/|\,\delta-1|$. 
Finally we use these results to establish self-adjointness criteria for degenerate elliptic diffusion operators.

\end{list}

\vfill



\vspace{0.5cm}

\noindent
\begin{tabular}{@{}cl@{\hspace{10mm}}cl}
$ {}^\dag\hspace{-5mm}$&   Mathematical Sciences Institute (CMA)    &  {} &{}\\
  &Australian National University& & {}\\
&Canberra, ACT 0200 && {} \\
  & Australia && {} \\
  &derek.robinson@anu.edu.au
 & &{}\\
\end{tabular}

\newpage

\setcounter{page}{1}

\section{Introduction}\label{S1}

Our intention is to analyse the $L_2$-version of the weighted Hardy inequality for functions with support  in a neighbourhood of the boundary $\Gamma$ of a domain $\Omega$ in $\Ri^d$
with special emphasis on evaluating the optimal constant in the inequality.
If   $d_\Gamma$ denotes the Euclidean distance to the boundary we assume that there are $\delta\geq0$ and $r_0>0$ and for each $r\in \langle0,r_0\rangle$
an $a>0$ such that
\begin{equation}
\|d_\Gamma^{\,\delta/2-1}\psi\|_2\leq a\,\|d_\Gamma^{\,\delta/2}\,(\nabla\psi)\|_2
\label{ehc1.1}
\end{equation}
for all $\psi\in C_c^1(\Gamma_{\!\!r})$ where $\Gamma_{\!\!r}=\{x\in\Omega:\,d_\Gamma(x)<r\}$.
The corresponding Hardy constant $a_\delta(\Gamma_{\!\!r})$  is then defined as the infimum of the 
$a$ for which (\ref{ehc1.1}) holds for all  $\psi\in C_c^1(\Gamma_{\!\!r})$.
Clearly  $a_\delta(\Gamma_{\!\!r})$ decreases as $r\to0$ and so we define the boundary constant $a_\delta(\Gamma)$ as the infimum over $r\in\langle0,r_0\rangle $  
of the $a_\delta(\Gamma_{\!\!r})$.
Although it is standard practise to examine the Hardy inequality (\ref{ehc1.1}) for all $\varphi\in C_c^1(\Omega)$ and to examine the corresponding optimal constant $a_\delta(\Omega)$ this approach has certain limitations since the inequality may fail even for very simple domains.
The general situation is clearly illustrated by the example of $\Omega=B(0\,;1)$, the unit ball.
Then (\ref{ehc1.1}) is valid for all $\varphi\in C_c^1(\Omega)$ if $\delta\in[0,1\rangle$ and $a_\delta(\Omega)=a_\delta(\Gamma)=2/(1-\delta)$.
If, however, $\delta>1$ then (\ref{ehc1.1}) is valid on $\Gamma_{\!\!r}$ for all $r\in\langle0,1\rangle$ but it fails on $\Omega$.
In this case one has $a_\delta(\Gamma)=2/(\delta-1)$.
The value $\delta=1$ is truly exceptional and the inequality also fails on all $\Gamma_{\!\!r}$. \cite{LamP}

The  boundary layer inequality (\ref{ehc1.1}) has occurred previously in several different contexts and it is increasingly clear that the value of the boundary Hardy constant has a special significance.
The inequality, with $\delta=0$,  was a principal ingredient in the paper of Brezis and Marcus \cite{BrM}, in which they analysed 
a one-parameter family of Hardy-type inequalities.
It also appeared with  small non-zero $\delta$  in Haj{\l}asz' paper on the pointwise Hardy inequality \cite{Haj1}
 and the subsequent work  of Filippas, Mazy'a and Tertikas \cite{FMT1} on Hardy-Sobolev inequalities.
The boundary constant also enters implicitly in  Nenciu and Nenciu's  self-adjointness criteria for Schr{\"o}dinger operators on domains \cite{NeN1} as 
explicitly established  by Ward \cite{War} \cite{War1}.
Inequality (\ref{ehc1.1}) also  arises naturally in the theory of elliptic operators  whose  coefficients have a boundary degeneracy of order $d_\Gamma^{\,\delta}$.
Then there is a unique critical value $\delta_c\in \langle 0,2]$, determined by the boundary constant $a_\delta(\Gamma)$, for which the operator   is self-adjoint  for all $\delta>\delta_c$ (see \cite{Rob16}, Theorem~5.2).
Another compelling reason to examine the inequality is that for  $\delta\in[0,2\rangle$ it is equivalent to a weighted version of the  weak Hardy  inequality introduced in 
\cite{Dav17} to analyse the local properties of the Hardy constant.
This equivalence will  be  demonstrated in Section~\ref{S2} for general $\Omega$  together with an extension of Davies' locality results to the weighted situation.

In Section~\ref{S3} we establish that the boundary inequality (\ref{ehc1.1}) has a universal trait. 
It is valid for all uniform domains, in the sense of Martio and Sarvas \cite{MarS}, with an Ahlfors regular boundary 
and for all $\delta\geq 0$ with the one exception $\delta=2-(d-d_{\!H})$ where $d_{\!H}$ denotes the Hausdorff dimension of $\Gamma$.
Moreover, the boundary constant $a_\delta(\Gamma)$ then satisfies the lower  bound
\[
a_\delta(\Gamma)\geq 2/|\,\delta-2 +d-d_{\!H}|
\;.
\]
The existence result is to a large extent a corollary of the results of Koskela and  Lehrb{\"a}ck  \cite{KL1} \cite{Leh3}  in combination with those of Haj{\l}asz \cite{Haj1}
 on the pointwise Hardy inequality.
The distinction between $\delta<2-(d-d_{\!H})$ and  $\delta>2-(d-d_{\!H})$ is an illustration of the dichotomy of Koskela and Zhong \cite{KZ} for the unweighted Hardy inequality. 
It is usually described as the difference between a thick and a thin boundary but in the theory of  degenerate elliptic operators it corresponds to small degeneracy and large degeneracy
at the boundary. 

In Sections~\ref{S4}, \ref{S5} and \ref{S6}  we pass from the general analysis  to the problem of evaluation of the boundary constant $a_\delta(\Gamma)$.
We successively discuss $C^{1,1}$-domains, convex domains and the complements of convex domains.
These cases  are  of particular   interest as they appear  to delineate the borderline between smooth and rough boundaries.
 The final picture is exactly the same in the first two cases.
 The boundary Hardy  inequality  (\ref{ehc1.1}) is satisfied for all $\delta\geq0$ with $\delta\neq1$ and the boundary constant $a_\delta(\Gamma)=2/|\,\delta-1|$.
 These results are known in part but there are several cases which have not apparently been analysed before.
 For example the conclusion
 for $C^{1,1}$-domains  appears to be new although it  is in part an extension of some known  results on $C^2$-domains.
 For example, Davies proved  for bounded  $C^2$-domains   that the optimal constant in the unweighted weak Hardy  inequality  is, with our convention,  
  equal to $2$ (see \cite{Dav17}, Theorems~2.3 and 2.5).
 Brezis and Marcus \cite{BrM}, Lemma~1.2,  subsequently established the same conclusion for the constant in the unweighted boundary inequality.
 Moreover, for general $C^2$-domains and $\delta>1$ the identification of $a_\delta(\Gamma)$ was established in \cite{Rob15}, Proposition~2.9.
 In Section~\ref{S5} we examine convex domains.
This is now a well studied situation for $\delta\in[0,1\rangle$.
 Matskewich and Sobolevskii \cite{MaSo}  proved  the optimal constant for the  
unweighted Hardy inequality on a convex $C^1$-domain is  equal to $2$.
 Davies (\cite{Dav17}, Theorem~2.6) also gave a straightforward proof   that  the optimal constant  in the weak version of the  inequality  is greater or equal to $2$ if $\Gamma$ has at least one point of convexity.
These works were followed by the proof of   Marcus, Mizel and Pinchover \cite{MMP}  that the constant is actually equal to $2$ for general convex domains and $\delta=0$.
The identification $a_\delta(\Gamma)=2/(1-\delta\,)$ for $\delta\in[0,1\rangle$ then follows immediately.
The definitive result for the small $\delta$ case is given by Avkhadiev \cite{Avk1}  together with references to earlier partial results.
It appears, however, that there are no previous results for convex domains and $\delta>1$.
In Section~\ref{S6} we assume $\Omega$ is the complement of the closure of a convex domain. 
Then the situation is more complicated.
If $\delta>1$ then the situation is similar to the above (see Theorems~1.1 and 4.1 in \cite{Rob13}). 
But for $\delta\in[0,1\rangle$ seemingly anomalous behaviour occurs even for the complements of convex polytopes, e.g. one can have $a_0(\Gamma)>2$.
This anomaly occurs if the boundary has sharply extruding edges.

Finally  in Section~\ref{S7} we use these results to derive self-adjointness criteria for the family of elliptic diffusion operators analysed   in \cite{Rob15} and \cite{Rob16}.
In particular we deduce   that $\delta>3/2$ is a sufficient condition for self-adjointness for $C^{1,1}$-domains and  the interior and exterior of convex domains.
It is also established that $\delta\geq 3/2$ is  a necessary condition.

\section{Hardy constants}\label{S2}

In this section we analyse  local geometric structure of the optimal  Hardy constant $a_\delta(\Gamma)$ corresponding to the boundary inequality  (\ref{ehc1.1}).
The starting point is a version of Davies' weak Hardy inequality \cite{Dav17} which  states that there are $b,c>0$ such that 
\begin{equation}
\|d_\Gamma^{\,\delta/2-1}\psi\|_2\leq b\,\|d_\Gamma^{\,\delta/2}\,(\nabla\psi)\|_2+c\,\|\psi\|_2
\label{ehc2.1}
\end{equation}
for all $\psi\in C_c^1(\Omega)$. 
This formulation differs from that of Davies insofar it is expressed in terms of the norms and not the squares of the norms.
But that is for convenience. 
The two formulations are equivalent for all subsequent purposes.

The weak Hardy constant $b_\delta(\Omega)$ corresponding to (\ref{ehc2.1}) is  defined as the infimum of the $b>0$ for which  there is a $c>0$ such that (\ref{ehc2.1}) is valid.
It is the square root of Davies' constant.
One can also define constants $b_\delta(\Gamma_{\!\!r})$  and $b_\delta(\Gamma)$ in analogy with the  $a_\delta(\Gamma_{\!\!r})$  and $a_\delta(\Gamma)$ by restriction to functions in $C_c^1(\Gamma_{\!\!r})$.
Moreover, if (\ref{ehc1.1}) extends to all of $C_c^1(\Omega)$ then one can also introduce the optimal constant  $a_\delta(\Omega)$.

Although the inequalities (\ref{ehc1.1}) nor (\ref{ehc2.1}) are ostensibly  distinct they are in fact equivalent for small $\delta$.
In particular the $\delta=0$ version of (\ref{ehc2.1}), the only case considered by Davies, is equivalent to the boundary 
inequality (\ref{ehc1.1}) with $\delta=0$.

 \begin{prop}\label{phc2.1}
Assume $\delta\in[0,2\rangle$.
Then the  boundary Hardy inequality $(\ref{ehc1.1})$ is valid if and only if the weak Hardy inequality $(\ref{ehc2.1})$ is valid.
Moreover, if the inequalities are satisfied then  $a_\delta(\Gamma)=b_\delta(\Gamma)=b_\delta(\Omega)$. 
Finally $a_\delta(\Gamma\cap U)=b_\delta(\Gamma\cap U)$ for each bounded open subset  $U\subset \Ri^d$.
\end{prop}

\noindent\proof\ (\ref{ehc1.1})$\,\Rightarrow \,$(\ref{ehc2.1})$\;$
Fix $\eta\in C^1(0,\infty)$ with $0\leq\eta\leq 1$, $\eta(s)=1$ if $s<1/2$, $\eta(s)=0$ if $s>1$ and $|\eta^{\,\prime}|\leq 3$.
Then define $\xi=\eta\circ(r^{-1}d_\Gamma)$. 
It follows that $0\leq \xi\leq 1$, $\xi=1$ on $\Gamma_{\!\!r/2}$, $\xi=0$ on $\Omega_r=\{x\in\Omega:d_\Gamma(x)>r\}$
and  $|\nabla\xi|\leq 3/r$.
Therefore since $1-\xi=0$ on $\Gamma_{\!\!r/2}$ one has
\begin{eqnarray*}
\|d_\Gamma^{\,\delta/2-1}\varphi\|_2&\leq& \|d_\Gamma^{\,\delta/2-1}\xi\varphi\|_2+\|d_\Gamma^{\,\delta/2-1}(1-\xi)\varphi\|_2\\[5pt]
&\leq& \|d_\Gamma^{\,\delta/2-1}\xi\varphi\|_2+(r/2)^{\delta/2-1}\|\varphi\|_2\leq \|d_\Gamma^{\,\delta/2-1}\xi\varphi\|_2+2\,r^{\delta/2-1}\|\varphi\|_2
\end{eqnarray*}
where the second step uses $\delta\in[0,2\rangle$.
Now, however, $\psi=\xi\varphi\in C_c^1(\Gamma_{\!\!r})$ and 
\[
\|d_\Gamma^{\,\delta/2-1}\xi\varphi\|_2\leq a\, \|d_\Gamma^{\,\delta/2}\,(\nabla(\xi\varphi))\|_2
\leq a\, \|d_\Gamma^{\,\delta/2}\,(\nabla\!\varphi)\|_2+3\,a r^{\delta/2-1}\|\varphi\|_2
\]
by (\ref{ehc1.1}) where we have used the Leibniz rule,  the bounds on $\xi$ and $|\nabla\xi|$ and  $\supp|\nabla\xi|\subseteq \Gamma_{\!\!r}$.
Then by combination of these two estimates one concludes that
\begin{equation}
\|d_\Gamma^{\,\delta/2-1}\varphi\|_2\leq a\,\|d_\Gamma^{\,\delta/2}(\nabla\!\varphi)\|_2+(2+3a)\,r^{\delta/2-1}\|\varphi\|_2
\label{ehc2.3}
\end{equation}
for all $\varphi\in C_c^1(\Omega)$, i.e.\  (\ref{ehc2.1})  is valid with $b=a$ and $c=(2+3a)\,r^{\delta/2-1}$.

\smallskip

\noindent (\ref{ehc2.1})$\,\Rightarrow \,$(\ref{ehc1.1})$\;$
If $r\in\langle0,r_0\rangle$ then (\ref{ehc2.1}) is automatically  valid for all $\psi\in C_c^1(\Gamma_{\!\!r})$.
But since $\delta<2$ one has  $r^{1-\delta/2}d_\Gamma^{\,\delta/2-1}\geq 1$ on $\Gamma_{\!\!r}$. 
Therefore
\[
\|d_\Gamma^{\,\delta/2-1}\psi\|_2\leq b\,\|d_\Gamma^{\,\delta/2}\,(\nabla\psi)\|_2+c\,r^{1-\delta/2}\|d_\Gamma^{\,\delta/2-1}\psi\|_2
\]
for all $\psi\in C_c^1(\Gamma_{\!\!r})$.
Then for each $\varepsilon\in\langle0,1\rangle$ one may choose $r$ sufficiently small that $c\,r^{1-\delta/2}<\varepsilon$.
In particular $(1-c\,r^{1-\delta/2})>1-\varepsilon>0$.
Therefore by rearrangement
\[
\|d_\Gamma^{\,\delta/2-1}\psi\|_2\leq(1-\varepsilon)^{-1}\,b\,\|d_\Gamma^{\,\delta/2}\,(\nabla\psi)\|_2
\]
for all $\psi\in C_c^1(\Gamma_{\!\!r})$.
Thus (\ref{ehc1.1}) is valid  on $\Gamma_{\!\!r}$ with $a=(1-\varepsilon)^{-1}\,b$.

\smallskip

Next consider the various Hardy constants.
It follows by definition that $b_\delta(\Gamma_{\!\!r})\leq a_\delta(\Gamma_{\!\!r})$ and $b_\delta(\Gamma_{\!\!r})\leq b_\delta(\Omega)$.
Therefore by taking infima over $r$ one obtains $b_\delta(\Gamma)\leq a_\delta(\Gamma)$ and  $b_\delta(\Gamma)\leq  b_\delta(\Omega)$.
(These bounds are valid for all $\delta\geq0$.)
Hence  to deduce that $a_\delta(\Gamma)=b_\delta(\Gamma)$ it suffices to prove that  $a_\delta(\Gamma)\leq b_\delta(\Gamma)$.

First for each  $\varepsilon\in\langle0,1\rangle$ one may choose  $r\in\langle0,r_0\rangle$ and $ c_{r,\varepsilon}>0$ such that 
\begin{equation}
\|d_\Gamma^{\,\delta/2-1}\psi\|_2\leq (1+\varepsilon)\,b_\delta(\Gamma)\,\|d_\Gamma^{\,\delta/2}\,(\nabla\psi)\|_2+c_{r,\varepsilon}\,\|\psi\|_2
\label{ehc2.4}
\end{equation}
for all $\psi\in C_c^1(\Gamma_{\!\!r})$.
In particular this is valid for all $\psi\in C_c^1(\Gamma_{\!\!s})$ with $s\in\langle0,r\rangle$.
But choosing $s$ such that $c_{r,\varepsilon}\,s^{1-\delta/2}<\varepsilon$ one deduces, as in 
 the proof of (\ref{ehc2.1})$\,\Rightarrow \,$(\ref{ehc1.1}),  that  (\ref{ehc1.1}) is valid on $\Gamma_{\!\!s}$
with $a=(1-\varepsilon)^{-1} (1+\varepsilon)\,b_\delta(\Gamma)$.
Therefore $a_\delta(\Gamma)\leq (1-\varepsilon)^{-1} (1+\varepsilon)\,b_\delta(\Gamma)$.
Hence in the limit $\varepsilon\to0$ one deduces that $a_\delta(\Gamma)\leq b_\delta(\Gamma)$.
Thus $a_\delta(\Gamma)=b_\delta(\Gamma)$.

The proof that $a_\delta(\Gamma\cap U)=b_\delta(\Gamma\cap U)$ is almost identical. 
The inequality $b_\delta(\Gamma\cap U)\leq a_\delta(\Gamma\cap U)$ follows  by definition and the converse follows by repeating the last argument
with $b_\delta(\Gamma)$ replaced by $b_\delta(\Gamma\cap U)$ and the $\psi$ restricted to $C_c^1(\Gamma_{\!\!s}\cap U)$.

Finally to deduce that $b_\delta(\Gamma)=b_\delta(\Omega)$  it suffices to prove that $b_\delta(\Omega)\leq b_\delta(\Gamma)$.
One can again assume (\ref{ehc2.4}) is valid for all $\psi\in C_c^1(\Gamma_{\!\!r})$.
Then  for each $\varphi\in C_c^1(\Omega)$ set $\varphi=\xi\varphi+(1-\xi)\varphi$ with $\xi$  the function introduced in the proof of (\ref{ehc1.1})$\,\Rightarrow\, $(\ref{ehc2.1}).
Then by a slight variation of that proof, using (\ref{ehc2.4}) with $\psi=\xi\varphi$, one establishes that   there is a $c_{r,\varepsilon}^{\,\prime}>0$ such that
\[
\|d_\Gamma^{\,\delta/2-1}\psi\|_2\leq (1+\varepsilon)\,b_\delta(\Gamma)\|d_\Gamma^{\,\delta/2}\,(\nabla\psi)\|_2+c^{\,\prime}_{r,\varepsilon}\,\|\psi\|_2
\]
for all $\psi\in C_c^1(\Omega)$.
Therefore $b_\delta(\Omega)\leq (1+\varepsilon)\,b_\delta(\Gamma)$ and in the limit $\varepsilon\to0$ one has $b_\delta(\Omega)\leq b_\delta(\Gamma)$.
Thus $b_\delta(\Omega)= b_\delta(\Gamma)$.
\hfill$\Box$

\bigskip

Next we turn to the characterization of $a_\delta(\Gamma)$ by the local constants $a_\delta(\Gamma\cap U)$.
It follows by definition that $a_\delta(\Gamma)\geq a_\delta(\Gamma\cap U)$.
Therefore
\begin{equation}
a_\delta(\Gamma)\geq \sup_{j\in \Ni}\,a_\delta(\Gamma\cap U_j)
\label{ehc2.5}
\end{equation}
where $\cu=\{U_j\}_{j\in\Ni}$ is a cover of $\Gamma$ by bounded open subsets $U_j$ of $\Ri^d$.
But, following Davies \cite{Dav17}, we establish  a converse inequality  by use of a suitably chosen partition of unity.

Let  $\cp=\{\psi_j\}_{j\in\Ni}$ be a locally finite partition of unity   by $C_c^1$-functions $\psi_j$ subordinate to the open cover $\cu$.
Thus $0\leq \psi_j\leq1$ and each  $\psi_j$ has support in a member of $\cu$ 
(see \cite{Rud2}, Theorem~6.20, or \cite{Spi}, Theorem~3.11).
Then  each $x\in \Ri^d$ has an open neighbourhood in which all but a finite number of the $\psi_j$ are zero,  $\sum_{j\in\Ni}\psi_j(x)=1$
and  the derivatives satisfy $\sum_{j\in\Ni}|\nabla\psi_j|<\infty$.
Set $\chi_j=\psi_j\,(\sum_{k\in\Ni}\psi_k^{\,2})^{-1/2}$.
It follows that $0\leq \chi_j\leq1$,  $\supp\chi_j= \supp\psi_j$ and the number of $\chi_j$ which are non-zero at any point $x$ is uniformly bounded.
But now one has $\sum_{j\in\Ni}\chi_j^{\,2}=1$ and  $\sum_{j\in\Ni}|\nabla\chi_j|^2<\infty$. (See  \cite{ReR}, pages 222 and 295, for a similar construction.)
 
\begin{prop}\label{phc2.2}
Assume $\delta\in[0,2\rangle$ and the boundary Hardy inequality $(\ref{ehc1.1})$  is satisfied. Then 
\[
a_\delta(\Gamma)=\sup_{j\in\Ni} a_\delta(\Gamma\cap U_j)
\;.
\]
\end{prop}
\proof\
First we establish the corresponding statement with $\Gamma$ replaced by $\Gamma\cap V$ where $V$ is a bounded open set.
Since $a_\delta(\Gamma\cap V)\geq \sup_{j\in\Ni} a_\delta(\Gamma\cap V\cap  U_j)$  it remains to prove the converse inequality.

Secondly,  fix $a_j$ such that (\ref{ehc1.1}) is valid with $a$ replaced by $a_j$ for all $\psi\in C_c^1(\Gamma_{\!\!r}\cap V\cap U_j)$.
Let $\Lambda=\{j\in\Ni: V\cap U_j\neq\emptyset\}$.
Then $\Lambda$  is finite.
Hence one may assume that the $a_j$ are uniformly bounded for all ${j\in\Lambda}$  and in fact for each $\varepsilon>0$ one may choose the $a_j$ such $a_j\leq (1+\varepsilon)\,a_\delta(\Gamma_{\!\!r}\cap V\cap U_j)$.
Now using the above partition of unity
\begin{eqnarray*}
\|d_\Gamma^{\,\delta/2-1}\psi\|_2^2&=&\sum_{j\in\Lambda}\|d_\Gamma^{\,\delta/2-1}\chi_j\psi\|_2^2
\leq \sum_{j\in\Lambda}a_j^{\,2}\,\|d_\Gamma^{\,\delta/2}\,(\nabla(\chi_j\psi))\|_2^2\\[5pt]
&\leq&(\max_{j\in\Lambda}a_j)^2\sum_{j\in\Lambda}\|d_\Gamma^{\,\delta/2}\,(\nabla(\chi_j\psi))\|_2^2
\end{eqnarray*}
for all $\psi\in C_c^1(\Gamma_{\!\!r}\cap V)$.
Thus if $a$ is the constant for which  (\ref{ehc1.1}) is valid for all $\psi\in C_c^1(\Gamma_{\!\!r}\cap V)$ then
one may assume that
\[
a\leq  (1+\varepsilon)\,\max_{j\in\Lambda}\,a_\delta(\Gamma_{\!\!r}\cap V\cap V_j)
\]
where we have used  $\supp\chi_j\psi\subset \Gamma_{\!\!r}\cap V\cap U_j$.

Next  it follows from Leibniz' rule and the Cauchy--Schwarz inequality that for each $\varepsilon_1>0$ one has an estimate
\begin{eqnarray*}
\sum_{j\in\Lambda}\|d_\Gamma^{\,\delta/2}(\nabla(\chi_j\psi))\|_2^2&\leq& (1+\varepsilon_1)\sum_{j\in\Lambda}\|d_\Gamma^{\,\delta/2}\,\chi_j(\nabla\psi)\|_2^2
+(1+\varepsilon_1^{-1})\sum_{j\in\Lambda}\||\nabla\chi_j|\, d_\Gamma^{\,\delta/2}\psi\|_2^2\\[5pt]
&\leq & (1+\varepsilon_1)\|d_\Gamma^{\,\delta/2}\,(\nabla\psi)\|_2^2+(1+\varepsilon_1^{-1})\,|\nabla\chi|^2\,r^\delta\| \psi\|_2^2
\end{eqnarray*}
where $ |\nabla\chi|^2=\sup_{x\in V}\sum_{j\in\Lambda}|(\nabla\chi_j)(x)|^2$.
This uses the properties of the partition of the unity.
Therefore by combining these estimates one obtains the weak Hardy inequality
\begin{eqnarray*}
\|d_\Gamma^{\,\delta/2-1}\psi\|_2^2&\leq& (1+\varepsilon_1)\, (\max_{j\in\Lambda}a_j)^2\,\|d_\Gamma^{\,\delta/2}\,(\nabla\psi)\|_2^2+(1+\varepsilon_1^{-1})\,|\nabla\chi|^2\,r^\delta \|\psi\|_2^2\\[5pt]
&\leq& \Big((1+\varepsilon_1)\, (\max_{j\in\Lambda}a_j)\,\|d_\Gamma^{\,\delta/2}\,(\nabla\psi)\|_2+(1+\varepsilon_1^{-1})\,|\nabla\chi|\,r^{\delta/2} \|\psi\|_2\Big)^2
\end{eqnarray*}
for all $\psi\in C_c^1(\Gamma_{\!\!r}\cap V)$.
Hence the weak Hardy constant $b_\delta(\Gamma_{\!\!r}\cap V)$ satisfies
\[
b_\delta(\Gamma_{\!\!r}\cap V)\leq (1+\varepsilon_1)\, (1+\varepsilon)\,\max_{j\in \Lambda}\,a_\delta(\Gamma_{\!\!r}\cap V\cap U_j)
\;.
\]
Taking the the infimum over $r$ followed by the  limits $\varepsilon_1\to0$ and $\varepsilon\to0$   one then concludes that 
\begin{equation}
b_\delta(\Gamma\cap V)\leq \sup_{j\in\Lambda}\,a_\delta(\Gamma\cap V \cap U_j)
\;.\label{ehc3.10}
\end{equation}
Since $\delta\in[0,2\rangle$ it follows by  Proposition~\ref{phc2.1} that $b_\delta(\Gamma\cap  V)=a_\delta(\Gamma\cap V)$.
Hence 
\[
a_\delta(\Gamma\cap V)\leq \sup_{j\in\Lambda}\,a_\delta(\Gamma\cap V\cap U_j)
\;.
\]
If $\Gamma$ is bounded, the proof of the proposition is complete
since for  sufficiently large $V$ one has $\Gamma\cap V=\Gamma$ and so $a_\delta(\Gamma)\leq \sup_{j\in\Lambda}\,a_\delta(\Gamma\cap U_j)$.

Finally the case of unbounded $\Omega$ is handled by another approximation argument.
Replace $V$ by a family of concentric open balls $B_n$ of radius $n$
and replace $\Lambda$ by $\Lambda_n$.
Then by the preceding
\[
a_\delta(\Gamma\cap B_n)\leq  \sup_{j\in\Lambda_n}\,a_\delta(\Gamma\cap B_n \cap U_j)
\leq \sup_{j\in\Ni}\,a_\delta(\Gamma\cap B_n\cap U_j)\leq a_\delta(\Gamma)
\;.
\]
Now we consider the limit as $n\to\infty$ which, by monotonicity,  coincides with the supremum over $n$.

Let $\chi_n\in [0,1]$ be an approximation to the identity on $\Ri^d$ with the property that $\chi_n=1$ on $B_n$, $\chi_n=0$ on the complement of $B_{2n}$ and $|\nabla\chi_n|\leq c\,n^{-1}$.
Then for each $\varepsilon>0$ one can assume that
\[
\|d_\Gamma^{\,\delta/2-1}(\chi_n\varphi)\|_2\leq (1+\varepsilon)\,a_\delta(\Gamma\cap  B_{2n})\,\|d_\Gamma^{\,\delta/2}\,\nabla(\chi_n\varphi)\|_2
\]
for all $\varphi\in C_c^1(\Gamma_{\!\!r})$.
But $a_\delta(\Gamma\cap  B_{2n})$ increases monotonically with $n$ so 
\[
\|d_\Gamma^{\,\delta/2-1}(\chi_n\varphi)\|_2\leq (1+\varepsilon)\,\Big(\,\textstyle{\sup_{m\in\Ni}}\,a_\delta(\Gamma\cap  B_{m})\Big)\,\|d_\Gamma^{\,\delta/2}\,\nabla(\chi_n\varphi)\|_2
\]
for all $\varphi\in C_c^1(\Gamma_{\!\!r})$.
Moreover,
\[
\|d_\Gamma^{\,\delta/2}\,\nabla(\chi_n\varphi)\|_2\leq \|d_\Gamma^{\,\delta/2}(\nabla\varphi)\|_2+r^{\delta/2}\|(\nabla\chi_n)\varphi\|_2
\;.
\]
Therefore
\[
\|d_\Gamma^{\,\delta/2-1}\chi_n\varphi\|_2\leq (1+\varepsilon)\,\Big(\,\textstyle{\sup_{m\in\Ni}}\,a_\delta(\Gamma\cap  B_{m})\Big)\,\|d_\Gamma^{\,\delta/2}\nabla(\varphi)\|_2+c\,n^{-1}\,r^{\delta/2}\|\varphi\|_2
\]
for all $\varphi\in C_c^1(\Gamma_{\!\!r})$.
Taking the limit $n\to\infty$ one deduces that  
\[
\|d_\Gamma^{\,\delta/2-1}\varphi\|_2\leq (1+\varepsilon)\,\Big(\,\textstyle{\sup_{m\in\Ni}}\,a_\delta(\Gamma\cap  B_{m})\Big)\,\|d_\Gamma^{\,\delta/2}\nabla(\varphi)\|_2
\]
for all $\varphi\in C_c^1(\Gamma_{\!\!r})$ and $\varepsilon>0$.
Hence
\[
a_\delta(\Gamma_{\!\!r})\leq \textstyle{\sup_{n\in\Ni}}\,a_\delta(\Gamma\cap  B_{n})
\;.
\]
Thus $a_\delta(\Gamma)\leq \sup_{n\in\Ni}a_\delta(\Gamma\cap  B_{n})$. 
As the converse inequality is evident from monotonicity one then has $a_\delta(\Gamma)= \sup_{n\in\Ni} a_\delta(\Gamma\cap B_n)$.

The foregoing argument can be repeated with $\Gamma$ replaced by $\Gamma\cap U_j$ to deduce that 
\[
a_\delta(\Gamma\cap U_j)=\sup_{n\in\Ni}a_\delta(\Gamma\cap U_j\cap B_n)
\;.
\]
Combining these observations with the conclusion of the previous paragraph one finds
\[
a_\delta(\Gamma)=\sup_{n\in\Ni}a_\delta(\Gamma\cap B_n)\leq \sup_{n\in\Ni}\sup_{j\in\Ni}\,a_\delta(\Gamma\cap B_n\cap U_j)
=\sup_{j\in\Ni}\,a_\delta(\Gamma\cap U_j)
\;.
\]
Since the converse inequality is given by (\ref{ehc2.5}) the proof of the proposition is complete.
\hfill$\Box$

\bigskip

After this preparation it is straightforward to establish the weighted version of Davies' local characterization of the Hardy constant.
\begin{thm}\label{thc2.1}
 Assume the weighted boundary Hardy inequality $(\ref{ehc1.1})$  is satisfied with  $\delta\in[0,2\rangle$.
Define the local Hardy constant $a_\delta(x)$ for  each $x\in\Gamma$ by
  \begin{equation}
  a_\delta(x)=\inf\{a_\delta(\Gamma\cap U): U\subset \Ri^d\,,\,U\ni x\}
  \label{ehc3.200}
  \end{equation}
  where the $U$ are bounded open subsets.
  It follows  $x\in\Gamma\mapsto a_\delta(x)$ is an upper semi-continuous function with values in $\langle0,a_\delta(\Gamma)]$
and 
  \[
  a_\delta(\Gamma\cap V)=\sup\{ a_\delta(x): x\in \Gamma\cap V\}
  \]
  for each bounded open subset $V$ of $\Ri^d$.
  Moreover, $a_\delta(\Gamma)=\sup\{ a_\delta(x): x\in \Gamma\}$.
  \end{thm}
  \proof\
  The theorem is a direct analogue of Lemma~2.2 and Theorem~2.3 in \cite{Dav17}.
  
First,   since $a_\delta(\Gamma\cap U)\wedge a_\delta(\Gamma\cap V)\geq  a_\delta(\Gamma\cap U\cap V)$
the value of $a_\delta(x)$ only depends on the geometry of an arbitrarily small neighbourhood of $x$.
 The function $x\in\Gamma\mapsto a_\delta(x)$ is strictly positive by definition and it follows from (\ref{ehc2.5}) that it is bounded by $a_\delta(\Gamma)$.
  The upper semi-continuity  follows directly from the  proof of Lemma~2.2 in \cite{Dav17}.

Secondly,  definition (\ref{ehc3.200}) gives  $ \sup\{ a_\delta(x): x\in \Gamma\cap U\}\leq a_\delta(\Gamma\cap U)$ for each bounded open  $U$.
But for each $\varepsilon>0$ and $x\in\Gamma\cap V$ there is a neighbourhood $U_x$ of $x$  such that $a_\delta(\Gamma\cap U_x)\leq \sup\{ a_\delta(y): y\in \Gamma\cap V\}+\varepsilon$.
  Then since $V$ is bounded there are a finite number of $x$ and corresponding $U_x$ such that these sets form an open cover of $\Gamma\cap V$.
  It then follows from Proposition~\ref{phc2.2} that $a_\delta(\Gamma\cap V)\leq \sup\{ a_\delta(y): y\in \Gamma\cap V\}+\varepsilon$.
  Taking the limit $\varepsilon\to0$ gives the second  statement of the theorem.
  The third, the identification of $a_\delta(\Gamma)$, follows by an approximation argument with an increasing family of $V$.
    \hfill$\Box$

\bigskip

It is unclear whether Theorem~\ref{thc2.1} extends to all $\delta\geq 2$ but if the weak weighted Hardy inequality  (\ref{ehc2.1}) is valid on $\Gamma_{\!\!r}$ for all small $r>0$
one obtains similar conclusions for the weak Hardy constant $b_\delta(\Gamma)$.
For example if $ b_\delta(x)=\inf\{b_\delta(\Gamma\cap U): U\subset \Ri^d\,,\,U\ni x\}$ then   $b_\delta(\Gamma\cap V)=\sup\{ b_\delta(x): x\in \Gamma\cap V\}$
  for each bounded open subset $V$ 
and in addition $b_\delta(\Gamma)=\sup\{ b_\delta(x): x\in \Gamma\}$.

\smallskip

Theorem~\ref{thc2.1} reduces the problem of calculating $a_\delta(\Gamma\cap V)$, or $a_\delta(\Gamma)$, to the calculation of $a_\delta(x)$, i.e.\ it reduces a global problem to a local one.
But the theorem depends on the assumption that $\delta\in[0,2\rangle$ and the weighted  boundary Hardy inequality is satisfied.
Since this is equivalent to Davies' weak Hardy inequality for this range of $\delta$ there appears to be no great gain.
The advantage, however, is that the boundary Hardy inequality is   valid under quite general circumstances
as demonstrated in the next section.

\section{Boundary Hardy inequality}\label{S3}

In this section we establish  a general theorem on the existence of the boundary Hardy inequality which is in part
pieced together from the literature on the Hardy inequality.
The principal difference with previous results is the extension to domains satisfying a locally uniform version of Ahlfors regularity.
The theorem has the advantage of allowing an estimation of a lower bound on the corresponding Hardy constant
$a_\delta(\Gamma)$, or the local constants $a_\delta(\Gamma\cap U)$, which is known to be optimal in many 
special cases.
The existence theorem depends on assuming that the domain $\Omega$ is uniform, in the sense of  Martio and Sarvas \cite{MarS}
and that the boundary $\Gamma$ satisfies an Ahlfors regularity property similar to that introduced in an earlier analysis of Markov uniqueness
of degenerate elliptic operators \cite{LR}.

First,   $\Omega$ is defined to be a uniform  domain if there is a $\sigma\geq1$ and  for each pair of points $x,y\in \Omega$  a rectifiable  curve 
$\gamma\colon[0,l]\to\Omega$, parametrized by arc length,  such that $\gamma(0)=x$, $\gamma(l)=y$  with arc length
$l(\gamma(x\,;y))\leq \sigma\,|x-y|$ and  $d_\Gamma(\gamma(t))\geq \sigma^{-1}\,(t\wedge(l-t))$ for all $t\in [0,l]$.
Uniform domains are  a special subclass of domains studied earlier by John \cite{John} in which the bound on the length of the curve $\gamma$ is omitted.
In fact these properties were initially only examined for bounded domains and the extension to 
unbounded domains was given subsequently by V{\"a}is{\"a}l{\"a} \cite{Vai1} (see also \cite{Leh3}, Section~4).

Secondly, there are various choices   for the 
Ahlfors regularity property of the  boundary~$\Gamma$.
The standard definition requires  that there exists a Borel measure $\mu$ on $\Gamma$,
a $C>1$ and an $s\in\langle 0,d\,]$ such that 
\[
C^{-1}\,r^s \leq \mu(\Gamma\cap {\overline B}(x\,;r)) \leq C\, r^s
\]
for each $x\in \Gamma$ and all $r\in\langle0,\diam(\Gamma)]$.
(As usual $B(x\,;r)$ is the open ball centred at $x$ with radius $r$ and ${\overline B}(x\,;r)$ its closure.)
A full analysis of this condition, expressed in a slightly different but equivalent form, can be found in \cite{KLV}, Chapter~7.
One immediate consequence of the condition is that $s=d_{\!H}$, the Hausdorff dimension of $\Gamma$,  and  $\mu$ is equivalent to the Hausdorff measure $\ch^s$ on $\Gamma$.
This regularity requirement is, however, a rather stringent condition if $\Gamma$ is unbounded and we will adopt the weaker locally uniform version  used in \cite{LR}.
Let $A_z=\Gamma\cap {\overline B}(z\,;R)$ with $z\in \Gamma$ and $R>0$.
Then $\Gamma$ is defined to be locally uniformly Ahlfors regular, or  uniformly Ahlfors $s$-regular, 
if there exists a Borel measure $\mu$ on $\Gamma$,
a $C>1$, an $R>0$  and an $s\in\langle 0,d\,]$ such that 
\begin{equation}
C^{-1}\,r^s \leq \mu(A_z\cap {\overline B}(x\,;r)) \leq C\, r^s
\label{ehrc3.1}
\end{equation}
for all $x\in A_z$, $r\in\langle0,2R\,]$ 
 and all $z\in\Gamma$.
The value of $C$ depends on $R$ but we have simplified notation by suppressing the parameter $R$.
Nevertheless,  if (\ref{ehrc3.1}) is satisfied for  all $x\in \Gamma\cap {\overline B}(z\,;R)$ and  $r\in\langle0,2R\,]$ 
 then it is satisfied for $x\in \Gamma\cap {\overline B}(z\,;S)$,  $r\in\langle0,2S\,]$ and $S<R$ with the same
 $C$ and $s$.
 Although the condition is weaker than the standard definition it still implies
 that $\mu$ and $\ch^s$  are  locally equivalent on $\Gamma$ and 
$s=d_{\!H}$.

\smallskip

The following theorem gives two conditions for the existence of the boundary Hardy inequality on all uniform domains with a locally uniform  Ahlfors regular boundary.
It concords with the dichotomy established by Koskela and Zhong \cite{KZ} for the unweighted Hardy inequality on a general domain.
The Koskela-Zhong conclusion  was later extended to  weighted Hardy inequalities by Lehrb{\"a}ck \cite{Leh3}.
Although the theorem is only stated on $L_2(\Omega)$ it has a direct analogue on $L_p(\Omega)$ for all $p\in\langle1,\infty\rangle$.

\begin{thm}\label{thc3.1}
Assume $\Omega$ is a uniform domain with a locally uniform Ahlfors regular boundary $\Gamma$.
Further assume that {\bf either } $\delta\in[0, 2-(d-d_{\!H})\rangle$ {\bf or} $\delta> 2-(d-d_{\!H})$.

Then there is an  $r_0>0$ and for each $r\in\langle0,r_0\rangle $ an $a>0$ such that
\begin{equation}
\|d_\Gamma^{\,\delta/2-1}\psi\|_2\leq a\,\|d_\Gamma^{\,\delta/2}\,(\nabla\psi)\|_2
\label{ehc3.1}
\end{equation}
for all $\psi\in C_c^1(\Gamma_{\!\!r})$.
\end{thm}
\proof\
The small $\delta$ case of the theorem is a version of the original result of Haj{\l}asz \cite{Haj1} who derived a boundary Hardy inequality from a  pointwise Hardy inequality.
His argument was restricted to the unweighted case but has subsequently been extended to the general weighted case (see \cite{KL1}).
In the weighted case  the  pointwise inequality states that there is a $C_1>0$ such that 
\begin{equation}
|\psi(x)|\leq C_1\,d_\Gamma(x)^{\,1-\delta/2}\,(M_{2d_\Gamma(x),q}(d_\Gamma^{\,\delta/2}\,|\nabla\psi|))(x)
\label{ehc3.2}
\end{equation}
for all $x\in\Gamma_{\!\!r}$ and $\psi\in C_c^1(\Gamma_{\!\!r})$ where $q\in \langle 1,2\rangle$.
Here 
\[
(M_{R,q}\,\psi)(x)=\sup_{r\in\langle0,R\rangle}\Big(|B(x\,;r)|^{-1}\int_{B(x\,;r)}dy\,|\psi(y)|^q\Big)^{1/q}
\]
is the local maximal operator.
Once one  establishes (\ref{ehc3.2}) the boundary Hardy inequality (\ref{ehc3.1}) follows  directly from the Hardy-Littlewood-Wiener maximal function theorem.
Specifically
\begin{equation}
\|d_\Gamma^{\,\delta/2-1}\psi\|_2\leq C_1\,\| M_{2d_\Gamma,q}(d_\Gamma^{\,\delta/2}\,|\nabla\psi|)\|_2\leq 
C_1\,C_2\,\,\|d_\Gamma^{\,\delta/2}\,(\nabla\psi)\|_2
\label{ehc3.20}
\end{equation}
where $C_2>0$ is the numerical constant in the $L_2$-Hardy-Littlewood-Wiener  theorem.
This depends only on $d$.
Thus the weighted Hardy inequality is valid on $L_2(\Gamma_{\!\!r})$.
Alternatively  if one can establish  (\ref{ehc3.2})  for all $\psi\in C_c^1(A_{z,r})$ where  $A_{z,r}=\Gamma_{\!\!r}\cap  {\overline B}(z\,;R)$
 then the weighted Hardy inequality  is valid on $L_2(A_{z,r})$.
Therefore the  strategy is to establish the pointwise inequality for all $A_{z,r}$ and to combine this information with a covering argument to obtain the weighted
 boundary Hardy inequality on $L_2(\Gamma_{\!\!r})$.
Thus the proof of the small $\delta$ statement in the theorem is now  in two steps, first the proof of the pointwise inequalities (\ref{ehc3.2})  in a suitably uniform fashion and 
secondly the extension of the local inequality (\ref{ehc3.20}) by a  covering argument.
The first step is a modification of the reasoning in \cite{KL1} and the second uses a partition of unity as in Section~\ref{S2}.

The proof of the pointwise inequality under  the current assumptions follows from a slight variation of the arguments in \cite{KL1}. In fact the principal
conclusion in the latter reference is valid under much more general hypotheses and does not require uniformity of the domain and does not involve
Ahlfors regularity. It relies, however,  on an estimate for  the  size of the boundary measured in terms of the Hausdorff content  $\ch^s_\infty$ of the boundary together with a  property of `visibility' of the boundary.
The required estimate on the Hausdorff content is given by the following.

\begin{lemma}\label{lrhc3.1}
It follows  from the locally uniform Ahlfors regularity hypothesis $(\ref{ehrc3.1})$ that there is a $C_{\!s}>1$
\begin{equation}
\ch^s_\infty(A_z\cap {\overline B}(x\,;r))\geq C_{\!s}^{\,-1} \,r^s
\label{ehc3.3}
\end{equation}
for all $x\in A_z$ and all $r\in\langle0,2R\,]$.
\end{lemma}
\proof\
Let $\{B(x_j\,;r_{\!j})\}_{j\in\Ni}$ be a cover of $A_z\cap B(x\,;r)$.
One may assume that there is a $y_j\in (A_z\cap B(x\,;r))\cap B(x_j\,;r_{\!j})$ for each $j\in \Ni$.
Then $B(x_j\,;r_{\!j})\subseteq B(y_j\,;2 r_{\!j})$ and so there is a second cover $\{B(y_j\,;t_{\!j})\}_{j\in\Ni}$  of $A_z\cap B(x\,;r)$ 
with the centres $y_j\in  A_z\cap B(x\,;r)$ and the radii $t_{\!j}=2r_{\!j}$.

If  $t_{\!j}>2R$ for some $j$ then
\[
\sum_{j\in\Ni}t_{\!j}^{\,s}>(2R)^s\geq r^s
\;.
\]
Consequently (\ref{ehc3.3}) is valid with $C_s=1$.
If, however, $t_{\!j}\leq 2R$  for all $j\in\Ni$ then   it follows  from the  regularity assumption that 
\[
C^{\,-1}r^s\leq \mu(A_z\cap B(x\,;r))\leq \sum_{j\in \Ni}\mu(B(y_j\,;t_{\!j})\leq C\,\sum_{j\in \Ni} t_{\!j}^{\,s}\
=2^sC\,\sum_{j\in\Ni}r_{\!\!j}^{\,s}
\;.
\]
where $C$ is the constant in the locally uniform Ahlfors regularity bounds (\ref{ehrc3.1}).
Thus taking the infimum over all such covers one 
deduces that 
\[
\ch^s(A_z\cap {\overline B}(x\,;r))\geq \ch^s_\infty(A_z\cap {\overline B}(x\,;r))\geq C_{\!s}^{\,-1}r^s
\;.
\]
Consequently (\ref{ehc3.3}) is valid with  $C_{\!s}=2^s\,C^2$.
\hfill$\Box$

\bigskip

If $\delta<2-(d-d_{\!H})$ the pointwise Hardy inequality (\ref{ehc3.2}) on  $C_c^1(A_{z,r})$  follows   for all small $r>0$  uniformly in $z$ from a slight modification of Proposition~4.3 and Theorem~5.1 in \cite{KL1}.
The proof is somewhat simplified  because of the assumption that the domain is uniform.
The new element of importance is  the uniformity in the position parameter $z$ occurring in the Ahlfors regularity property.
The pointwise inequality involves a constant $C_1$ whose value is determined by the constant in the Hausdorff content estimate (\ref{ehc3.3}), the inner boundary 
condition in the terminology of \cite{KL1}.
The proof of Lemma~\ref{lrhc3.1} demonstrates that this is dependent solely on the constant $C$ in the local  Ahlfors regularity condition (\ref{ehrc3.1}) and the regularity 
parameter $s$,  which is equal to the Hausdorff dimension $d_{\!H}$  of the boundary.
Then the constant in  the  resulting boundary Hardy inequality (\ref{ehc3.20}) also involves the constant $C_2$ in the Hardy--Littlewood--Wiener theorem
which only depends on the dimension $d$.
Therefore the $C_1$ in the   boundary Hardy inequality is determined entirely by the Ahlfors constant $C$, the dimension $d$ of the space and the Hausdorff dimension $d_{\!H}$ of
the domain.
It is uniform in $z$ for all the sections $A_{z,r}$ of the boundary.
Consequently (\ref{ehc3.20})  is uniform for all the subspaces $L_2(A_{z,r})$.

The slight modifications of the proof of \cite{KL1}, other than the simplifications which occur due to the uniformity of the domain, involve the locality.
Condition~4.1 of \cite{KL1} involves the Hausdorff content of a section of the boundary, the visual section near $x$,  and that is replaced by the 
bounds  (\ref{ehc3.3}) on the sections $A_z\cap {\overline B}(x\,;r)$
and $r$ is identified with  $d_\Gamma(x)$. 
The restriction on the range of $r$ in the Ahlfors property (\ref{ehrc3.1}) and the Hausdorff content estimate (\ref{ehc3.3}) 
then restricts the conclusion of Theorem~5.1 in \cite{KL1} to small values of $d_\Gamma(x)$, i.e.\ the pointwise Hardy inequality is restricted to a boundary layer.
We will not delve into the details of the modifications as the proof of Theorem~5.1 is rather detailed but the end result is that the pointwise Hardy inequality
(\ref{ehc3.2}) is valid for all  $x\in A_{z,r}$ and $\psi\in C_c^1(A_{z,r})$ with $A_{z,r}=\Gamma_{\!\!r}\cap {\overline B}(z\,;R)$ and the boundary Hardy inequality is valid 
 on the  subspaces $L_2(A_{z,r})$  uniformly for $z\in \Gamma$ and $r\in\langle0,r_0]$ for some $r_0>0$.

The second step in the proof of the small $\delta$ case is the extension of the local Hardy inequality
(\ref{ehc3.20}) from $L_2(A_{z,r})$ to $L_2(\Gamma_{\!\!r})$.
First we construct a cover of $\Gamma_{\!\!r}$ by a countable family of balls $B(x_j\,;r_{\!j})$ with centres $x_j\in \Gamma$ and uniformly bounded radii $r_{\!j}$.
This is straightforward.
If $\{B(x_j\,;r_{\!j})\}_{j\in\Ni}$  is a cover of $\Ri^d$ with the $r_{\!j}\in\langle 0, \rho\rangle $ one eliminates all balls with $B(x_j\,;r_{\!j})\cap \Gamma_{\!\!r}=\emptyset$.
Then in the remaining balls one chooses $y_j\in B(x_j\,;r_{\!j})\cap \Gamma_{\!r}$ and introduces the balls $B(y_j\,;2r_{\!j}+r_0)$.
These new balls constitute a cover of $\Gamma_{\!\!r}$ for all $r\in\langle0,r_0\rangle$ and $B(y_j\,;2r_{\!j}+r_0)\cap \Gamma\neq\emptyset$.
Finally choose $z_j\in B(y_j\,;2r_{\!j}+r_0)\cap \Gamma$ and set  $\rho_j=2r_{\!j}+r_0$.
Then the  family of balls $ B(z_j\,;\rho_j)$ cover $\Gamma_{\!\!r}$, the centres $z_j$ are in $\Gamma$ and the radii satisfy $\rho_j<R$ where $R=2\rho+r_0$.
In particular $\Gamma_{\!\!r}\cap B(z_j\,;\rho_j)\subset A_{z_j,r}$ where $A_{z_j,r}=\Gamma_{\!\!r}\cap {\overline B}(z_j\,;R)$.

Next by refining the cover, if necessary, one may assume that  each point $x\in \Gamma_{\!\!r}$ is contained in at most a finite number $N$ of balls in  the cover.
After the refinement one may choose a partition of unity $\{\psi_j\}_{j\in\Ni}$ with $C_c^1$-functions $\psi_j$ 
as in Section~\ref{S2}.
Thus $0\leq\psi_j\leq1$, $\supp\psi_j\in B(x_j\,;r_{\!j})$, $ \sum_{j\in\Ni} \psi_j(x)=1$ and  $ \sum_{j\in\Ni} |(\nabla\psi_j)(x)|<\infty$  for all $x\in\Gamma_{\!\!r}$.
Finally one sets $\chi_j=\psi_j\,(\sum_{k\in\Ni}\psi_k^{\,2})^{-1/2}$.
It follows that $0\leq \chi_j\leq1$,  $\supp\chi_j= \supp\psi_j$ and all but $N$  of the $\chi_j$ are zero at any point $x$.
But now one has $\sum_{j\in\Ni}\chi_j(x)^{\,2}=1$ and  $\sum_{j\in\Ni}|(\nabla\chi_j)(x)|^2<\infty$ for all $x\in\Gamma_{\!\!r}$.
Each sum is finite for each $x\in\Gamma_{\!\!r}$ with a maximum of $N$-terms.
Therefore if $\varphi\in C_c^1(\Gamma_{\!\!r})$ then $\chi_j\varphi\in C_c^1(A_{z_j;r})$.
Consequently $\psi=\chi_j\varphi$ satisfies the pointwise Hardy inequality (\ref{ehc3.2}) and consequently the bounded Hardy inequality (\ref{ehc3.20})
on $L_2(A_{z_j,r})$.
Specifically one has a family of inequalities
\begin{equation}
\|d_\Gamma^{\,\delta/2-1}(\chi_j\varphi)\|_2^2\leq 
C^2_3\,\|d_\Gamma^{\,\delta/2}\,(\nabla(\chi_j\varphi))\|_2^2
\label{ehc3.4}
\end{equation}
where the constant $C_3$ depends on the Ahlfors constant $C$, the dimensions $d$ and $d_{\!H}$ and the localization length $R$ but not on $z_j$.
Thus $C_3$ is independent of $j$.
Hence summing the inequalities and using the properties of the partition of unity one finds
\[
\|d_\Gamma^{\,\delta/2-1}\varphi\|_2^2=\sum_{j\in\Ni}\|d_\Gamma^{\,\delta/2-1}(\chi_j\varphi)\|_2^2\leq 
C^2_3\,\sum_{j\in\Ni}\|d_\Gamma^{\,\delta/2}\,(\nabla(\chi_j\varphi))\|_2^2
\;.
\]
Then by the Leibniz' rule and the Schwarz inequality one deduces that for each $\varepsilon>0$
\[
\|d_\Gamma^{\,\delta/2-1}\varphi\|_2^2\leq 
C^2_3\,(1+\varepsilon)\,\|d_\Gamma^{\,\delta/2}\,(\nabla\varphi)\|_2^2
+C^2_3\,(1+\varepsilon^{-1})\sum_{j\in\Ni}\|d_\Gamma^{\,\delta/2}\,(\nabla\chi_j)\varphi\|_2^2
\;.
\]
But $\sup_{j\in\Ni}\|\nabla\chi_j\|_\infty\leq K<\infty$ and $\supp|\nabla\chi_j|\subset \supp|\chi_j|$.
Since $0\leq\chi_j\leq 1$ and  at most $N$ of the $\chi_j$ have overlapping support it follows that 
\[
\sum_{j\in\Ni}\|d_\Gamma^{\,\delta/2}\,(\nabla\chi_j)\varphi\|_2^2\leq  K^2N^2\,\|d_\Gamma^{\,\delta/2}\,\varphi\|_2^2
\leq K^2N^2r^2\,\|d_\Gamma^{\,\delta/2-1}\,\varphi\|_2^2
\]
where the last step uses $\supp\varphi\subset \Gamma_{\!\!r}$.
Therefore
\[
\|d_\Gamma^{\,\delta/2-1}\varphi\|_2^2\leq C_3^2\Big((1+\varepsilon)\,\|d_\Gamma^{\,\delta/2}\,(\nabla\varphi)\|_2
+(1+\varepsilon^{-1})KNr\,\|d_\Gamma^{\,\delta/2-1}\,\varphi\|_2\Big)^2
\;.
\]
Then by combination and rearrangement of these latter two estimates one obtains the boundary Hardy inequality
\[
\|d_\Gamma^{\,\delta/2-1}\varphi\|_2\leq C_r\,(1+\varepsilon)\,\|d_\Gamma^{\,\delta/2}\,(\nabla\varphi)\|_2
\]
on $\Gamma_{\!\!r}$
with $C_r=C_3\,\left(1-C_3\,(1+\varepsilon^{-1})KNr\right)^{-1}$ whenever $C_3\,(1+\varepsilon^{-1})KNr<1$.
Thus $a_\delta(\Gamma)=(1+\varepsilon)\inf_{r}C_r=(1+\varepsilon)C_3$ for all $\varepsilon>0$.
Hence $a_\delta(\Gamma)=C_3$.

\smallskip

This completes the proof of the first statement of the theorem. 
The second statement is essentially a corollary of  Theorem~1.3 in \cite{Leh3}.

First the locally uniform version (\ref{ehrc3.1}) of the Ahlfors regularity suffices to prove that the Aikawa dimension of the boundary 
is equal to the Hausdorff dimension by the proof  of Lemma~2.1 in \cite{Leh3}.
Therefore Theorem~1.3 in \cite{Leh3} states that if $\Omega$ is an unbounded John domain and $\delta>2-(d-d_{\!H})$
then the weighted Hardy inequality is valid on $\Omega$.
But each uniform  domain is a John domain so this proves the second statement of the theorem for unbounded $\Omega$.
If, however, $\Omega$ is bounded a slight modification of the proof of \cite{Leh3}, Theorem~1.3, as explained in \cite{Rob16}, Section~6,
establishes the boundary Hardy inequality.

It should be noted that the proof of the boundary inequality for large $\delta$ is quite different to the proof for small $\delta$.
For small $\delta$ the proof is based on the existence of the pointwise Hardy inequality.
But it follows from Theorem~1.1 in \cite{Leh6} that the condition $\delta<2-(d-d_{\!H})$ is necessary for the pointwise Hardy  inequality.
So an analogous argument for large $\delta$ is not possible.
\hfill$\Box$

\bigskip

Theorem~\ref{thc3.1} demonstrates that the boundary Hardy inequality (\ref{ehc3.1}), the $L_2(\Gamma_{\!\!r})$-inequality, has a universal character
not shared by the standard $L_2(\Omega)$-inequality.
The boundary inequality is valid in  many situations for which the full inequality fails.
But the arguments of Koskela and Lehrb{\"a}ck \cite{KL1} \cite{Leh3} \cite{Leh6} also give conditions which ensure that the $L_2(\Omega)$-inequality is valid.
For example, if $\delta>2-(d-d_{\!H})$ then, as mentioned above, the standard Hardy inequality on $\Omega$ is valid for all unbounded uniform domains with a locally uniform Ahlfors regular boundary.
If, however, $\delta<2-(d-d_{\!H})$ then the full inequality follows for these domains  if the pointwise inequality (\ref{ehc3.2}) is valid for all $x\in\Omega$.
Alternatively, it follows if the visual part of the boundary (see \cite{KL1}, Condition~4.1) is  the complete boundary.
Then one has the Hausdorff content condition 
\[
\ch^s_\infty(\Gamma\cap {\overline B}(x\,;r))\geq C_s^{-1}r^s
\]
for all $x\in\Gamma$ and $r>0$.

\smallskip

There is an interesting domain distinction between the two cases of the theorem.
\begin{remarkn}\label{rhc3.1}
If $\delta\in[0,2-(d-d_{\!H})\rangle$ then the boundary Hardy inequality (\ref{ehc3.1}) follows from the pointwise inequality (\ref{ehc3.2}).
Thus if the pointwise inequality is valid for all $\psi\in C^1_c(\Gamma_{\!\!r})$ and $s<r$ then the boundary inequality is satisfied on $\Gamma_{\!\!s}$ for functions which
do not necessarily vanish on the inner boundary $\Gamma_{\!\!s}^{(i)}=\{x\in\Omega:d_\Gamma(x)=s\}$.
Thus after closure the inequality (\ref{ehc3.1})  is valid for all weakly differentiable $\psi\in L_2(\Gamma_{\!\!s})$ with $d_\Gamma^{\,\delta/2}|\nabla\psi|\in L_2(\Gamma_{\!\!s})$
satisfying a Dirichlet boundary condition on the outer boundary $\Gamma$ but with no restriction on the inner boundary $\Gamma_{\!\!s}^{(i)}$.
On the other hand if $\delta>2-(d-d_{\!H})$ then the arguments of \cite{Leh3} establish (\ref{ehc3.1}) for all $\psi\in C_c^1(\Gamma_{\!\!s})$.
Then it follows by closure that (\ref{ehc3.1}) extends to all $\psi\in L_2(\Gamma_{\!\!s})$ with $d_\Gamma^{\,\delta/2}|\nabla\psi|\in L_2(\Gamma_{\!\!s})$
satisfying a Dirichlet boundary condition on the inner boundary  $\Gamma_{\!\!s}^{(i)}$ but with no restriction on the outer  boundary $\Gamma$.
This is a consequence of the large degeneracy  on $\Gamma$ which ensures that the weighted capacity of $\Gamma$ is zero (see \cite{LR}, Section~3).
\end{remarkn}

Although Theorem~\ref{thc3.1} establishes the  boundary Hardy inequality for all $\delta\geq0$, with the exception of the value $\delta=2-(d-d_{\!H})$, it does not 
give any precise  information on the  boundary Hardy constant $a_\delta(\Gamma)$.
But one can derive a lower bound on the constant by a variation of the arguments of \cite{Rob16}.
In fact one does not need to assume the  Ahlfors regularity property (\ref{ehrc3.1})  but only a property of the volume growth near the boundary which is a consequence of the regularity.
If $A$ is a bounded subset of $\Gamma$ and $A_r=\{ x\in\overline{\Omega}:d(x\,;A)<r\}$ then it follows from (\ref{ehrc3.1}) that  there are $\kappa, \kappa^{\,\prime}>0$
such that  
\begin{equation}
\kappa^{\,\prime} r^{(d-d_{\!H})}\leq |A_r|\leq \kappa \,r^{(d-d_{\!H})}
\label{eqhc3.1}
\end{equation}
for all small $r$
where $|A_r|$ denotes the Lebesgue measure of $A_r$ and $d_{\!H}$ is the Hausdorff dimension of $ A$.
In addition a mild form of the uniformity property is also required  if  $d_{\!H}\in[d-1,d\,\rangle$ (see \cite{LR}, Section~2, for details).

Now  for a general domain $\Omega$ we define
$x\in \Gamma$ to be an Ahlfors point if there is a bounded open $A=\Gamma\cap B(x\,;s)$  satisfying (\ref{eqhc3.1})
for all small $s>0$.

The following result is a version of Proposition~6.4 in \cite{Rob16} which was stated with the restriction $\delta>2-(d-d_{\!H})$.
It is, however, valid for all $\delta\geq 0$ with the exception of $\delta=2-(d-d_{\!H})$.

\begin{prop}\label{prhc3.1}
Assume  the weighted boundary Hardy inequality $(\ref{ehc1.1})$ is satisfied on $\Gamma_{\!\!r}\cap U$ where $U$ is a bounded open
 subset of $\Ri^d$ containing  an Ahlfors point  $x$ of  the boundary $\Gamma$. 
Set $a_\delta(x)=\inf\{a_\delta(\Gamma\cap V): U\supseteq V\ni x\}$ then
\[
 a_\delta(x)\geq 2/|(d-d_{\!H})+\delta-2\,|
\]
for  all $\delta\geq 0$ with  $\delta\neq 2-(d-d_{\!H})$.
\end{prop}
\proof\
The proof   is similar to the proof of the earlier result in \cite{Rob16}.

Set $\alpha=\beta+\delta-2$ where $\beta =d-d_{\!H}$.
Then replace $\psi\in C_c^1(\Gamma_{\!\!r}\cap U)$ in the Hardy inequality (\ref{ehc1.1}) by $d_\Gamma^{\,-\alpha/2}\psi$.
It follows that 
\begin{eqnarray*}
\|d_\Gamma^{\,-\beta/2}\psi\|_2&=&\|d_\Gamma^{\,\delta/2-1}(d_\Gamma^{\,-\alpha/2}\psi)\|_2\leq
a\,\|d_\Gamma^{\,\delta/2}\,\nabla(d_\Gamma^{\,-\alpha/2}\psi)\|_2\\[5pt]
&\leq &a\,\|d_\Gamma^{\,\delta/2}\,(\nabla d_\Gamma^{\,-\alpha/2})\psi\|_2+a\,\|d_\Gamma^{\,\delta/2}d_\Gamma^{\,-\alpha/2}\,(\nabla \psi)\|_2\\[5pt]
&\leq&a\,|(\beta+\delta-2)/2\,|\,\|d_\Gamma^{\,-\beta/2}\psi\|_2+a\,\|d_\Gamma^{\,1-\beta/2}(\nabla\psi)\|_2
\end{eqnarray*}
for all $\psi\in C_c^1(\Gamma_{\!\!r}\cap U)$.
Therefore
\begin{equation}
1\leq a\,|\beta+\delta-2\,|/2+ a\,\|d_\Gamma^{\,1-\beta/2}\,(\nabla\psi)\|_2/\|d_\Gamma^{\,-\beta/2}\psi\|_2
\label{eqhc3.3}
\end{equation}
for all $\psi\in C_c^1(\Gamma_{\!\!r}\cap U)$.

Next  one constructs a sequence of $\psi_n\in C_c^1(\Gamma_{\!\!r}\cap U)$ such that   the numerator in the last term 
of (\ref{eqhc3.3}) is bounded uniformly in $n$ but the denominator diverges as $n\to\infty$.
(Since $\beta$ is independent of $\delta$ the value of $\delta$ plays no role in this latter argument.)
 Once this is achieved one immediately deduces that 
 $1\leq a\,|\,\beta+\delta-2\,|/2$.
Therefore  optimizing over the choice of $a$ one obtains
\[
1\leq a_\delta(\Gamma_{\!\!r}\cap U)\,|\,\beta+\delta-2\,|/2
\]
for all $r\in\langle0,r_0\rangle$.
Finally taking the infimum over subsets $V$ of $U$ containing $ x$ one deduces  that $1\leq a_\delta(x)\,|\,\beta+\delta-2\,|/2$.

\smallskip

 The   $\psi_n$ are constructed as  in \cite{Rob16}.
First define $\xi_n\in W^{1,\infty}(0,\infty)$ by
$\xi_n(t)=0$ if $t<1/n$, $\xi_n(t)=1$ if $t>1$,  $\xi_n(t)=\log(nt)/\log n$ if $1/n\leq t\leq 1$ and $\xi_n(t)=1$ if $t>1$.
Thus  $\xi_n(t)\geq 1/2$ if $t\in [n^{-1/2},1]$.
Then fix a decreasing function  $\chi\in C^1(0,1)$ with $\chi(t)=1$ if $t\in[0,1/2]$, $\chi(1)=0$ and $|\chi'|\leq 3$.
Finally define $\psi_n=(\xi_n\circ(r^{-1}d_\Gamma))(\chi\circ d_A )$  where $A=\Gamma\cap U$ and $d_A(x)=d(x\,;A)$.
It follows  that  $0\leq\psi_n\leq 1$, $\supp\psi_n\subset A_1\cap \Gamma_{\!\!r}$ and  the $\psi_n$ converge pointwise to $\chi\circ d_A$ as $n\to\infty$.
Further $(\xi_n\circ(r^{-1}d_\Gamma))\geq 1/2$ if $r^{-1}d_\Gamma\in[n^{-1/2},1]$ and $(\chi\circ d_A)\geq1$ if $d_A\leq1/2$.
In addition $d_A\geq d_\Gamma$.

Therefore if $r\leq 1$
\[
\|d_\Gamma^{\,-\beta/2}\psi_n\|_2=\int d_\Gamma^{\,-\beta}\,|\psi_n|^2\geq (1/4)\int d_A^{\,-\beta}\,\one_{D_{r, n}}
\]
where $D_{r, n}=\{x\in\Omega:r/n^{1/2}\leq d_A(x)\leq r/2\,,\,r/n^{1/2}\leq d_\Gamma(x)\leq r\}$ and we have 
used the inclusion $\{x\in \Omega: d_A(x)<1/2\}\supset \{x\in \Omega: r/n^{1/2}\leq d_A(x)<r/2\}$.
Then since $d_A^{\,-\beta}=r^{-\beta}(1+\beta\int^1_{d_A(x)/r}t^{-(\beta+1)})$ one has 
\begin{eqnarray*}
\int d_\Gamma^{\,-\beta}\,|\psi_n|^2&\geq &(1/4)\,r^{-\beta}\int \one_{D_{r,n}}\Big(1+\beta\int^1_{d_A(x)/r}t^{-(\beta+1)}\Big)\\[5pt]
&\geq &(1/4)\,r^{-\beta}|D_{r,n}| +(1/4)\,\beta\int^1_{ n^{-1/2}}dt\,t^{-1}((rt)^{-\beta}|D_{r,t,n}|)
\end{eqnarray*}
where $D_{r,t, n}=\{x\in\Omega:r/n^{1/2}\leq d_A(x)\leq tr/2\,,\,r/n^{1/2}\leq d_\Gamma(x)\leq r\}$.
But since $U\ni x$ and $x$ is an Ahlfors point it follows from (\ref{eqhc3.1}) that $\lim_{n\to\infty}(r^{-\beta}|D_{r,n}|)=(r^{-\beta}|A_{r/2}|)\geq 2^{-\beta}\kappa'$ and 
$\lim_{n\to\infty}((rt)^{-\beta}|D_{r,t,n}|)=((rt)^{-\beta}|A_{rt/2}|)\geq 2^{-\beta}\kappa'$ where the bound  is
uniform for $t\in\langle0,1\rangle$.
Since the integral of $t^{-1}$ is divergent at the origin one concludes that $\|d_\Gamma^{\,-\beta/2}\psi_n\|_2\to \infty$ as $n\to\infty$.
This is the first step in handling the last term in (\ref{eqhc3.3}).

The second step is to estimate $\|d_\Gamma^{\,1-\beta/2}(\nabla\psi_n)\|_2^2=\int d_\Gamma^{\,2-\beta}|\nabla\psi_n|^2$.
First observe that
\[
|\nabla\psi_n|^2\leq 2r^{-2}|(\xi'_n\circ d_\Gamma)|^2\, |(\chi\circ d_A )|^2
+2r^{-2}|(\xi_n\circ(r^{-1}d_\Gamma))|^2\,|(\chi'\circ d_A )|^2
\]
by the Leibniz rule and the Cauchy-Schwarz inequality.
Denote the integrals involving the first and second terms on the right hand side by $I_1$ and $I_2$, respectively.

First assume $\beta\leq2$.
Then since $\supp \psi_n\subset \Gamma_{\!\!r}$, $|\xi_n|\leq 1$  and $|\chi'|\leq 3$ one has 
\begin{eqnarray*}
I_2\leq 18r^{-2}\int d_\Gamma^{\,2-\beta}\one_{\{x:\,0< d_\Gamma(x)\leq r\}} \,\one_{\{x:\,0< d_{A}(x)\leq 1\}}
\leq 18r^{-\beta}\int \one_{\{x:\,0< d_B(x)\leq r\}}
\end{eqnarray*}
where $B\supseteq A$ is a slight enlargement of $A$. 
Then $I_2\leq 18\kappa$, uniformly for all small $r$, by the bounds (\ref{eqhc3.1}) applied to $B$.

Alternatively, if $\beta>2$ one deduces that 
\[
I_2\leq 18r^{-2}\int d_\Gamma^{\,2-\beta}\, \one_{D'_{r,n}}\leq 18\Big(r^{-\beta}\,|D'_{r,n}|+(\beta-2)\int^1_{n^{-1}}dt\,t\,((rt)^{-\beta}\,|D'_{r,t,n}|\Big)
\]
where $D'_{r,t,n}=\{x\in\Omega:rn^{-1}\leq d_\Gamma(x)\leq rt\,,\,1/2\leq d_A(x)\leq 1\}$ and $D'_{r,n}=D'_{r,1,n}$.
First  $\lim_{n\to\infty}(r^{-\beta}|D'_{r,n}|)\leq \kappa\,|A_1|$.
Secondly $D'_{r,t,n}\subset\{x\in\Omega:0\leq d_\Gamma(x)\leq rt\,,\,0< d_A(x)\leq 1\}$.
Therefore $D'_{r,t,n}\subset { B}_{rt}$ where $B$ is again a slight  enlargement of $A$.
Hence $(rt)^{-\beta}|D'_{r,t,n}|\leq (rt)^{-\beta}|B_{rt}|$ is bounded uniformly for all $n\geq 1$ and $t\leq 1$.
Consequently $I_2$ is uniformly bounded in $n$ uniformly for all small $r$.

Finally we have to estimate the integral $I_1$.
But 
\[
I_1=2\,r^{-2}\int d_\Gamma^{\,2-\beta}|(\xi'_n\circ(r^{-1}d_\Gamma))|^2\,|(\chi\circ d_A )|^2\leq 2(\log n)^{-2}\int d_\Gamma^{\,-\beta}\,\one_{D''_{r,n}}
\]
where one now has  $D''_{r,n}=\{x\in\Omega:rn^{-1}\leq d_\Gamma(x)\leq r\,,\,0< d_A(x)\leq 1\}$.
Arguing as above
\[
I_1\leq 2\,(\log n)^{-2}\Big(|D''_{r,n}|+\beta\int^1_{n^{-1}} dt\,t^{-1}((rt)^{-\beta}\,|D''_{r,t,n}|\Big)
\]
with $D''_{r,t,n}=\{x\in\Omega:rn^{-1}\leq d_\Gamma(x)\leq rt\,,\, 0< d_A(x)\leq 1\}$.
The  Ahlfors assumption (\ref{eqhc3.1}) implies  that $\sup_{n\geq1}(r^{-\beta}|D''_{r,n}|)\leq \kappa$ and  $\sup_{n\geq1}((rt)^{-\beta}|D''_{r,t,n}|)\leq \kappa$.
Since  the integral then  gives a factor $\log n$ one obtains a bound $I_1\leq b\, (\log n)^{-1}$ with $b>0$.
Hence one concludes that
\[
 \|d_\Gamma^{\,1-\beta/2}(\nabla\psi_n)\|_2=\int d_\Gamma^{\,2-\beta}|\nabla\psi_n|^2\leq a+b\,(\log n)^{-1}
 \]
 with $a,b>0$.
 This completes the proof of  Proposition~\ref{prhc3.1}.
 \hfill$\Box$

\begin{cor}\label{cohc3.1}
If  the weighted boundary Hardy inequality $(\ref{ehc1.1})$ is satisfied on $\Gamma_{\!\!r}$ with $\delta\geq 0$ and if there is an Ahlfors point  $x\in\Gamma$ then
\[
a_\delta(\Gamma)\geq 2/|(d-d_{\!H})+\delta-2\,|
\]
for all  $\delta\geq 0$ with the exception of $\delta=2-(d-d_{\!H})$.
\end{cor}
\proof\
It follows directly from the Hardy inequality that $a_\delta(\Gamma)\geq a_\delta(\Gamma\cap U)$ for all bounded open $U\subset \Ri^d$.
But then $a_\delta(\Gamma)\geq a_\delta(x)\geq 2/|(d-d_{\!H})+\delta-2\,|$ by Proposition~\ref{prhc3.1}.
\hfill $\Box$

\bigskip

Note that the conclusion of the corollary is valid for all uniform domains with boundaries satisfying the locally uniform Ahlfors regularity property
if $\delta\neq 2-(d-d_{\!H})$ by Theorem~\ref{thc3.1}.

\begin{remarkn}\label{remhc3.1} 
Corollary~\ref{cohc3.1} sheds some light on Davies' conjecture \cite{Dav17} on the Hardy constant.
With our definition of the weak constant Davies' conjectured `that if the weak unweighted Hardy inequality is valid
on the bounded domain $\Omega$ then $b_0(\Omega)\geq2$'.
But $b_0(\Omega)=a_0(\Gamma)$ by Proposition~\ref{phc2.1} so Corollary~\ref{cohc3.1} implies that $b_0(\Omega)\geq 2/|d-d_{\!H}-2|$.
Since $\Omega$ is bounded $d_{\!H}\geq d-1$, the topological dimension of $\Gamma$,  and then the bound gives $b_0(\Omega)\geq 2$ if and
only if $d_{\!H}=d-1$.

If $d_{\!H}>d-1$ it is possible that $b_0(\Omega)<2$.
A simple example is the bounded two-dimensional domain with boundary the von Koch curve constructed from an equilateral triangle.
Then $d=2$, the boundary is Ahlfors regular and $d_{\!H}=\log 4/\log3$ and Corollary~\ref{cohc3.1} asserts that
$a_0(\Gamma)=b_0(\Omega)\geq \log4/\log3\in\langle1,2\rangle$.
It is feasible that this bound is attained and Davies' conjecture is false.
\end{remarkn}

Although the lower bound of Corollary~\ref{cohc3.1} is not attained in general (see Example~6.9 in \cite{Rob16}) it is possible that it 
is attained if the boundary is a self-similar fractal.
In the next two sections we establish that the lower bound of the corollary has a matching upper bound for $C^{1,1}$-domains or general convex domains
and all $\delta\geq0$.
Thus in these cases the lower bound is attained.

\section{$C^{1,1}$-domains}\label{S4}

In this section we consider the weighted boundary Hardy inequality for $C^{1,1}$-domains and, more generally,
for domains  whose boundaries have a point with a $C^{1,1}$-neighbourhood.
The validity  of the boundary Hardy inequality  for  $C^{1,1}$-domains  can be deduced from  Theorem~\ref{thc3.1} with $d_{\!H}=d-1$.
Nevertheless  we give an independent proof which exploits the details of the $C^{1,1}$-property and leads to the identification of the optimal
Hardy constant as $2/|\,\delta-1|$.
Thus the bound of Corollary~\ref{cohc3.1} is attained.

The $C^{1,1}$-property is a   bound on the curvature of the boundary.
If $\Omega$ is bounded then it is a $C^{1,1}$-domain if and only if it satisfies both a
 uniform interior ball condition and a uniform exterior  ball condition.
Specifically the interior  condition requires that for each $x\in \Gamma$ there exists a  $y\in\Omega$ and a $u>0$ such that 
${\overline B}(y\,;u)\cap \Omega^{\,\rm c}=x$.
The exterior condition is  defined similarly by interchanging $\Omega$ and $\Omega^{\,\rm c}$.
It follows by definition that the principal curvatures at each point of the boundary are bounded uniformly by $u^{-1}$.
These conditions have important implications for the distance to the boundary $d_\Gamma$  both in the interior and the exterior of 
$\Omega$.

The two-sided ball condition  implies that each point $y\in \Gamma_{\!\!u}$, where $\Gamma_{\!\!u}$ is either an interior or exterior boundary layer,  has a unique nearest boundary  point $x=n(y)\in\Gamma$.
Consequently  
$(\nabla d_\Gamma)(y)=(y-x)/|y-x|$ (see, for example, \cite{BEL}, Theorem~2.2.7).
In particular $|(\nabla d_\Gamma)(y)|=1$.
Secondly, $d_\Gamma\in C^{1,1}(\Gamma_{\!\!u})$ and this implies that 
 the partial derivatives $\partial_jd_\Gamma$ of $d_\Gamma$ are locally weakly$^*$ differentiable.
Therefore $d_\Gamma\in W^{2,\infty}_{\rm loc}(\Gamma_{\!\!u})$.
This allows one to obtain estimates on  the Hessian $D^2d_\Gamma=(\,\partial_j\partial_k d_\Gamma\,)$ of $d_\Gamma$
analogous to those established originally for $C^{\,2}$-domains in \cite{GT}, Appendix~14.6.
In the $C^{\,2}$-case one has
\begin{equation}
|(\nabla^2d_\Gamma)(y)|\leq \Tr((|D^2d_\Gamma|)(y))\leq \gamma\,u^{-1}
\label{ehc4.1}
\end{equation}
for all $y\in\Gamma_{\!\!r}$ with $r<u/2$ where the value of $\gamma$ depends only on the dimension $d$.
(One can choose $\gamma=2\,(d-1)$ by the following lemma.)
The corresponding $C^{1,1}$-estimates are local and only require a local version of the $C^{1,1}$-property.

A point $x\in \Gamma$ is defined to be a $C^{1,1}$-point
 if there is a bounded open subset $U\subset \Ri^d$  such that $U\ni x$ and  the section $\Gamma\cap U$ of the boundary is the graph of a $C^{1,1}$-function, after a   suitable affine change of coordinates (see, for example, \cite{GT}, Section~6.2).
 Then the two-sided ball condition is satisfied locally, i.e.\ one can choose  $U\subset  \Ri^d$ and $ u>0$ such that
  $d_\Gamma$ is a $C^{1,1}$-function on $\Gamma_{\!\!u}\cap U$.

\begin{lemma}\label{lhc4.1}
 Assume $x$ is a $C^{1,1}$-point of $\Gamma$ so $d_\Gamma\in W^{2,\infty}_{\rm loc}(\Gamma_{\!\!u}\cap U)$.
 Then there is a $\gamma>0$ such that 
 \[
\Big| \sum^d_{j=1}(\partial_j \psi, \partial_jd_\Gamma)\Big|
\leq \gamma\,u^{-1}\,\|\psi\|_1
\]
for all $\psi\in C_c^1(\Gamma_{\!\!r}\cap U)$ with $r<u/2$.
One may choose $\gamma=2\,(d-1)$.
\end{lemma}
\proof\
The proof is essentially the same as the argument in the $C^2$-case.
One obtains by the calculations of \cite{GT} the estimates
\begin{eqnarray*}
\Big| \sum^d_{j=1}(\partial_j \psi, \partial_jd_\Gamma)\Big|&=&
\Big|\int_{\Gamma_{\!\!r}\cap U}dx\, \psi(x)(\nabla^2d_\Gamma)\,\Big|\\[-2pt]
&=&\Big|\int_{\Gamma_{\!\!r}\cap U}dx\, \psi(x)\sum^{d-1}_{j=1}\kappa_j(x)\,(1-d_\Gamma(x)\kappa_j(x))^{-1}\Big|\\[-2pt]
&\leq &\int_{\Gamma_{\!\!r}\cap U}dx\, |\psi(x)|\sum^{d-1}_{j=1}|\kappa_j(x)|\,(1-d_\Gamma(x)|\kappa_j(x)|)^{-1}
\end{eqnarray*}
for $\psi\in C_c^1(\Gamma_{\!\!r}\cap U)$ where $\kappa_j(x)$ are the principal curvatures of $\Gamma$ at the unique nearest point $n(x)\in \Gamma$ of $x\in \Gamma_{\!\!r}$.
But $|\kappa_j(x)|\leq u^{-1}$ and $d_\Gamma(x)\leq r<u/2$ so the statement of the lemma follows immediately with $\gamma=2\,(d-1)$.
\hfill$\Box$

\bigskip
The actual value of the dimension dependent constant $\gamma$ is not important. 
In fact reduction of $r$ leads to a reduction of $\gamma$ and $\gamma\to (d-1)$ as $r\to0$.
The estimate of the lemma nevertheless  allows the derivation of a local  version of the weighted Hardy inequality and calculation of  the optimal 
local constant at the $C^{1,1}$-point $x$.

\begin{thm}\label{trhc4.1}
Assume  the boundary $\Gamma$ of the domain $\Omega$ contains a $C^{1,1}$-point $x\in U$.
Then the weighted boundary Hardy inequality $(\ref{ehc1.1})$ is satisfied on $\Gamma_{\!\!r}\cap U$ for all small $r>0$.
Moreover, if  $a_\delta(x)=\inf\{a_\delta(\Gamma\cap U): U\ni x\}$ then
\[
a_\delta(x)=2/|\,\delta-1|
\]
for all $\delta\geq 0$ with  $\delta\neq 1$.
\end{thm}
\proof\
The proof of  the local version of the Hardy inequality is based on the identity
\begin{equation}
(\delta-1)\,(d_\Gamma^{\,\delta-2}\varphi^2)\,(\nabla d_\Gamma)^2=(\nabla d_\Gamma).( \nabla(d_\Gamma^{\,\delta-1}\varphi^2))
-d_\Gamma^{\,\delta-1}((\nabla d_\Gamma).(\nabla\varphi^2)
\;.
\label{ihc1}
\end{equation}
Then since $x$ is a $C^{1,1}$-point it follows that $|(\nabla d_\Gamma)(y)|=1$ for all $y\in \Gamma_{\!\!u}\cap U$.
Therefore
\[
|\,\delta-1|\,\|d_\Gamma^{\,\delta/2-1}\varphi\|_2^2\leq \Big|\int_{{\Gamma_{\!\!u}\cap U}}(\nabla d_\Gamma).( \nabla(d_\Gamma^{\,\delta-1}\varphi^2))\Big|
+2\,\Big|\int_{{\Gamma_{\!\!u}\cap U}} d_\Gamma^{\,\delta-1}\varphi\,((\nabla d_\Gamma).(\nabla\varphi))\Big|
\]
for  $\varphi\in C^1_c(\Gamma_{\!\!u}\cap U)$.
But by Lemma~\ref{lhc4.1} 
\[
\Big|\int_{{\Gamma_{\!\!r}\cap U}}(\nabla d_\Gamma).( \nabla(d_\Gamma^{\,\delta-1}\varphi^2)\Big|\leq 
\gamma\,u^{-1}\,\|d_\Gamma (d_\Gamma^{\,\delta-2}\varphi^2)\|_1= \gamma\,(r/u)\,\|d_\Gamma^{\,\delta/2-1}\varphi\|_2^2
\]
for all $\varphi\in C^1_c(\Gamma_{\!\!r}\cap U)$ with $r<u/2$ and $\gamma=2\,(d-1)$.
Therefore, substituting this estimate in the identity and rearranging, one obtains 
\begin{eqnarray*}
(|\,\delta-1|-\gamma\,(r/u))\|d_\Gamma^{\,\delta/2-1}\varphi\|_2^2&\leq & \int_{{\Gamma_{\!\!r}\cap U}}d_\Gamma^{\,\delta-1}\,|(\nabla d_\Gamma).(\nabla\varphi^2)|
\leq 2\,\|d_\Gamma^{\,\delta/2}\,(\nabla\varphi)\|_2
\,\|d_\Gamma^{\,\delta/2-1}\varphi\|_2
\;.
\end{eqnarray*}
Since $\delta\neq1$  one can choose $r>0$ such that $\gamma\,(r/u)<|\,\delta-1|$.
Hence
\[
(|\,\delta-1|-\gamma\,(r/u))\,\|d_\Gamma^{\,\delta/2-1}\varphi\|_2\leq 2\,\|d_\Gamma^{\,\delta/2}\,\nabla\varphi\|_2
\]
for all $\varphi\in C_c^1(\Gamma_{\!\!r}\cap U)$ where $ U\ni x$.
Thus the weighted Hardy inequality is valid on $\Gamma_{\!\!r}\cap U$  for all small $r$  and 
\[
a_\delta(\Gamma\cap U)\leq 2/|\,\delta-1|
\;.
\]
Hence  $a_\delta(x)\leq 2/|\,\delta-1|$.
Finally as $x$ is a $C^{1,1}$-point of $\Gamma$ it is also an Ahlfors point of $\Gamma$ with $d_{\!H}=d-1$.
Therefore $a_\delta(x)\geq 2/|\,\delta-1|$ by Proposition~\ref{prhc3.1}. Consequently  $a_\delta(x)= 2/|\,\delta-1|$.
\hfill$\Box$

\bigskip

The local statement of Theorem~\ref{trhc4.1}  can be extended to a similar result for general $C^{1,1}$-domains but one has
to be precise about the definition of the $C^{1,1}$-property  in the case of unbounded domains.
In the bounded case the standard definition (see \cite{GT}, Section~6.2) is by local diffeomorphisms.
But this definition is equivalent to the two-sided uniform ball condition (see \cite{Barb} or \cite{Dalp}).
Alternatively, the definition is equivalent to the signed distance function being  a $C^{1,1}$-function in a neighbourhood of the 
boundary (see \cite{DeZ}, Section~7.8).
Each of these definitions can be extended to unbounded domains  in an equi-continuous manner (see \cite{DeZ}, Chapter~2) 
which respects the various equivalences.
 In order to avoid  justification of these various equivalences
 we will define an unbounded 
 domain $\Omega$ to be a $C^{1,1}$-domain if  and only if the uniform interior and exterior 
ball conditions are satisfied.
This ensures that the curvature of $\Gamma$ is uniformly bounded and the foregoing local estimates are valid.

\begin{thm}\label{trhc4.2}
If $\Omega$ is a  $C^{1,1}$-domain  then the boundary Hardy inequality $(\ref{ehc1.1})$ is satisfied for all small $r$ and all $\delta\geq0$ with $\delta\neq 1$.

Moreover, the boundary Hardy constant is given by
$a_\delta(\Gamma)=2/|\,\delta-1|$.
\end{thm}
\proof\
Fix $\varphi\in C_c(\Gamma_{\!\!r})$.
Then let $\{\chi_j\}_{j\in\Ni}$ denote the partition of unity  introduced in Section~\ref{S2}.
It follows that $\chi_j\in C_c^1(\Gamma_{\!\!r}\cap U_j)$.
Thus the estimates established in the proof of Theorem~\ref{trhc4.1} can be applied to each $\chi_j\varphi$.
In particular one has
\begin{eqnarray*}
(|\,\delta-1|-\gamma\,(r/u))\|d_\Gamma^{\,\delta/2-1}(\chi_j\varphi)\|_2&\leq& 2\,\|d_\Gamma^{\,\delta/2}\,\nabla(\chi_j\varphi)\|_2
\end{eqnarray*}
for all $j\in\Ni$ and all $r>0$ satisfying $\gamma\,(r/u)<|\,\delta-1|$.
But using the Leibniz rule and proceeding as in the proof of Theorem~\ref{thc3.1} one deduces that for each $\varepsilon>0$ one has
\[
(|\,\delta-1|-\gamma\,(r/u))^2\,\|d_\Gamma^{\,\delta/2-1}(\chi_j\varphi)\|_2^2\leq 
4\,(1+\varepsilon)\,\|\chi_j\,d_\Gamma^{\,\delta/2}\,(\nabla\varphi)\|_2^2
+4\,(1+\varepsilon^{-1})\,\|(\nabla\chi_j)\,d_\Gamma^{\,\delta/2}\varphi\|_2^2
\;.
\]
Now there is a $K>0$ such that $\sup_{j\in\Ni}|\nabla\chi_j|\leq K$,  by the definition of the partition of unity.
In addition $\supp|\nabla\chi_j|\subseteq \supp \chi_j$.
Therefore letting $\one_j$ denote the characteristic function of $\supp\chi_j$  and using the basic property of the partition of unity one obtains
\[
(|\,\delta-1|-\gamma\,(r/u))^2\,\|d_\Gamma^{\,\delta/2-1}\varphi\|_2^2\leq 4\,(1+\varepsilon)\,\|d_\Gamma^{\,\delta/2}\,(\nabla\varphi)\|_2^2
+4\,(1+\varepsilon^{-1})\,K^2\sum_{j\in\Ni}\|\one_j\,d_\Gamma^{\,\delta/2}\varphi\|_2^2
\;.
\]
But each $x\in \Gamma_{\!\!r}$ is contained in at most a finite number $N$ of the sets $\supp\chi_j$.
Hence
\begin{eqnarray*}
(|\,\delta-1|-\gamma\,(r/u))^2\,\|d_\Gamma^{\,\delta/2-1}\varphi\|_2^2&\leq& 4\,(1+\varepsilon)\,\|d_\Gamma^{\,\delta/2}\,(\nabla\varphi)\|_2^2
+4\,(1+\varepsilon^{-1})\,K^2\,N\,\|d_\Gamma^{\,\delta/2}\varphi\|_2^2\\[5pt]
&\leq&4\,\left((1+\varepsilon)\|d_\Gamma^{\,\delta/2}\,(\nabla\varphi)\|_2+(1+\varepsilon^{-1})KN^{1/2}\|d_\Gamma^{\,\delta/2}\varphi\|_2\right)^2
\;.
\end{eqnarray*}
Therefore, setting $\gamma_1=2(1+\varepsilon^{-1})KN^{1/2}u$ and noting that
$\|d_\Gamma^{\,\delta/2}\varphi\|_2\leq r\,\|d_\Gamma^{\,\delta/2-1}\varphi\|_2$ 
one deduces that 
\[
(|\,\delta-1|-(\gamma+\gamma_1)\,(r/u))\,\|d_\Gamma^{\,\delta/2-1}\varphi\|_2\leq 2\,(1+\varepsilon)\,\|d_\Gamma^{\,\delta/2}\varphi\|_2
\]
for all $r>0$ such that $(\gamma+\gamma_1)\,(r/u)<|\,\delta-1|$.
Since this is valid for all $\varphi\in C_c^1(\Gamma_{\!\!r})$  it follows that the weighted boundary Hardy inequality is valid on the 
boundary layer  $\Gamma_{\!\!r}$.
Moreover,  $a_\delta(\Gamma_{\!\!r})\leq 2(1+\varepsilon)/(|\,\delta-1|-(\gamma+\gamma_1)\,(r/u))$.
Therefore $a_\delta(\Gamma)\leq 2(1+\varepsilon)/|\,\delta-1|$ for all $\varepsilon>0$ and  for all $\delta\geq 0$ with the exception of $\delta=1$.
Hence $a_\delta(\Gamma)\leq 2/|\,\delta-1|$. 

Finally as $x$ is a $C^{1,1}$-point of $\Gamma$ it is also an Ahlfors point of $\Gamma$ with $d_{\!H}=d-1$.
Therefore $a_\delta(x)\geq 2/|\,\delta-1|$ by Proposition~\ref{prhc3.1}. Consequently  $a_\delta(x)= 2/|\,\delta-1|$.
\hfill$\Box$

\bigskip

The theorem extends the conclusion obtained for $C^{\,2}$-domains  in \cite{Rob13}, Section~2.4.
The latter relied on the stronger estimate (\ref{ehc4.1}) and was restricted to the case $\delta>1$.
The $C^{1,1}$-case is of greater interest as it marks the borderline at which any argument relying
on the twice-differentiability of $d_\Gamma$ fails.
For example, if $\Omega$ is a $C^{1,\alpha}$-domain with $\alpha\in\langle0,1\rangle$ then $d_\Gamma$ is
at most a $C^{1,\alpha}$-function in a neighbourhood of the boundary.

\section{Convex domains}\label{S5}

In this section we consider the evaluation of the boundary Hardy constant for convex domains $\Omega$.
This is  well understood for the unweighted case.
Then the  Hardy inequality is  valid on all convex domains and, with our convention,  the Hardy constant $a_0(\Omega)=2$ (see  \cite{MMP}, Appendix~A).
Hence it  follows that
\[
\|d_\Gamma^{\,\delta/2-1}\varphi\|_2\leq 2\,\|\nabla(d_\Gamma^{\,\delta/2}\varphi)\|_2
\leq 2\,\|d_\Gamma^{\,\delta/2}(\nabla\varphi)\|_2+2\,(\delta/2)\,\|d_\Gamma^{\,\delta/2-1}\varphi\|_2
\]
for all $\varphi\in C_c^1(\Omega)$.
If $\delta<1$ one then deduces, by rearrangement,  that the weighted Hardy inequality is  valid on $\Omega$ with $a_\delta(\Omega)\leq 2/|\,\delta-1|$.
In particular $a_\delta(\Gamma_{\!\!r})\leq 2/|\,\delta-1|$ for all small $r>0$ and consequently $a_\delta(\Gamma)\leq 2/|\,\delta-1|$.
But since $d_{\!H}=d-1$ and  each point of $\Gamma$ is an Ahlfors point one has  $a_\delta(\Gamma)\geq 2/|\,\delta-1|$
by Proposition~\ref{prhc3.1}.
Therefore $a_\delta(\Gamma)=2/|\,\delta-1|$ for all $\delta\in[0,1\rangle$.
(More complete results  on $L_p$-inequalities can be found in  \cite{Avk1} together with references to various  earlier results
on convex domains.)

If $\delta>1$ then the situation is different and seems not to have been explored. 
First one cannot expect the weighted Hardy inequality to be valid on the whole domain.
In fact it fails on the unit ball $B(0\,;1)$ (see, for example, \cite{Rob16} Example~5.6).
Nevertheless each convex domain is a uniform domain with a $(d-1)$-Ahlfors regular boundary.
Therefore the weighted Hardy inequality is valid at least on a boundary layer $\Gamma_{\!\!r}$ for all $\delta>1$ by
Theorem~\ref{thc3.1}.
Hence it remains to calculate the boundary constant $a_\delta(\Gamma)$.
Since $a_\delta(\Gamma)\geq 2/|\,\delta-1|$,
by Proposition~\ref{prhc3.1},
a matching upper bound is required.
We will achieve this by  approximation of  $\Omega$ by an increasing family of convex subdomains.
The idea of  monotonic approximation by a  family of special convex subdomains has a long history.
For example, Hadwiger \cite{Had} and  Eggleston \cite{Egg} derive approximations  of this type with families of regular convex subsets or convex polytopes
and  Grisvard \cite{Gris}, Lemma~3.2.1.1, states without proof an approximation theorem in terms of $C^2$-subdomains.
Barb (see \cite{Barb}, Theorem~5.1.33) also constructs an interesting approximation in terms of $C^{1,1}$-domains.  
The first calculation  \cite{MMP} of the Hardy constant for  the unweighted Hardy inequality  on a general convex domain was  based on approximations with convex polytopes.
This reasoning was also used  by Brezis and Marcus \cite{BrM}, Section~5, in their extension of the Hardy inequality.
We continue the exploitation of this idea.

\begin{thm}\label{tc5}
If $\Omega$ is a convex domain then the boundary Hardy inequality $(\ref{ehc1.1})$ is satisfied for all small $r$ and all $\delta\geq0$ with $\delta\neq 1$.

Moreover, the boundary Hardy constant is given by $a_\delta(\Gamma)=2/|\,\delta-1|$.
\end{thm}
\proof\
First the case $\delta\in[0,1\rangle$ is handled by the results of  \cite{MMP}, Appendix~A, for $\delta=0$ and the argument in the first paragraph above.
Secondly, if $\delta>1$  the existence of the boundary Hardy inequality 
 follows from Theorem~\ref{thc3.1} and the lower bound on $a_\delta(\Gamma)$ from
Proposition~\ref{prhc3.1}.
It remains to establish the upper bound on $a_\delta(\Gamma)$ and this will be achieved by elaboration of the arguments of \cite{MMP} and \cite{BrM}, Section~6.
We  adopt much of  the notation of the latter references.

\smallskip

Let $S$ be a bounded convex polytope in $\Ri^d$ and $\Gamma_{\!\!1},\ldots ,\Gamma_{\!\!n}$ the open $(d-1)$-dimensional faces of $S$.
Thus the boundary  $\Gamma=\bigcup^{\,n}_{j=1}{\overline\Gamma}_{\!\!j}$.
Let $\pi_j$ denote the hyperplane containing $\Gamma_{\!\!j}$ and $G_j$ the half space containing $S$ such that $\partial G_j=\pi_j$.
Then $S=\bigcap^{n}_{j=1} G_j$.
If $x\in\Ri^d$ let $d_j(x)=d(x\,;\pi_j)$ and $n_j(x)$ the unique nearest point of $x$ in $\pi_j$.
Further let $S_j$ denote the open subsets  defined by
\[
S_j=\{x\in S: d_j(x)<d_k(x) \mbox{ for all } k\neq j\}
\;.
\]
It is established in \cite{BrM}, Section~5, that $d_\Gamma(x)=\min(d_1(x), \ldots d_n(x))$ and $d_\Gamma(x)=d_j(x)$ implies 
that $n_j(x)\in {\overline \Gamma}_{\!\!j}$.

The boundary layer ${ \Gamma}_{\!\!r}$ is,  for all small $r$,  a set sandwiched between   an  exterior boundary component   $\Gamma$ and an  interior component
 ${\widetilde\Gamma}_{\!\!r}=\{x\in\Omega: d_\Gamma(x)=r\}$.
Next consider the open subsets $\Gamma_{\!\!j,r}=\{x\in S_j: d_j(x)<r\}$  of the boundary layer $\Gamma_{\!\!r}$ associated with the $S_j$.
Then the boundary $\partial\Gamma_{\!\!j,r}$  of $\Gamma_{\!\!j,r}$  is the union of the closures of the exterior face $\Gamma_{\!\!j}$, the analogous interior face 
${\widetilde\Gamma}_{\!\!j}=\{x\in{\widetilde\Gamma}_{\!\!r}:d_j(x)=r\}$ and the interfaces $I_{\!j,k}$, $k\neq j$, where
\[
I_{\!j,k}=\{x\in \Gamma_{\!\!r}: d_j(x)=d_k(x)\}
\;.
\]
The first step in the remaining proof is the analysis of the boundary Hardy inequality corresponding to $S$.
  
  \begin{prop}\label{pcnotes1}
  If $S$ is a bounded convex polytope then   $a_\delta(\Gamma)=2/(\delta-1)$ for all $\delta>1$.
   \end{prop}
  \proof\
  First we consider  the halfspaces $G_j$ for orientation.
Then the problem is essentially  one-dimensional and can be resolved with the aid of the inequality
    \begin{equation}
\int^r_0 dt\, t^\delta f'(t)^2\geq ((\delta-1)/2)^2 \int^r_0dt\,t^{\delta-2}f(t)^2-((\delta-1)/2)\,r^{\delta-1}f(r)^2
\label{ehcl1.0}
\end{equation}
which is valid for all $\delta\geq 0$ and all  $f\in W^{1,2}(0,r)$ with $t\in\langle0,r\rangle\mapsto t^{-1}f(t)\in L_2(0,r)$.
 (The proof of this inequality is contained in the multi-dimensional calculation below.)
The inequality already indicates the difference between the cases $\delta<1$ and $\delta>1$.
If $\delta<1$ the boundary term on the right hand side is positive and can be discarded to give a standard Hardy inequality.
But if $\delta>1$ the boundary term is negative and has to be retained.
If, however, $f(r)=0$ it again plays no role.
This simplification occurs  for the half-spaces $G_j$.
Choosing coordinates such that $G_j=\{ (x', x_j) : x'\in \Ri^{d-1},\,x_j>0\}$  one  deduces from (\ref{ehcl1.0}) that 
\begin{eqnarray*}
\int_{G_j}dx\, d_j(x)^\delta\,|(\nabla\varphi)(x)|^2&\geq& \int_{\Ri^{d-1}}dx'\int^\infty_0 dx_j \,x_j^\delta\,|(\partial_j\varphi)(x',x_j)|^2\\
&=&((\delta-1)/2)^2 \int_{\Ri^{d-1}}dx'\int^\infty_0 dx_j \,x_j^{\delta-2}\,|\varphi(x',x_j)|^2\\
&=&((\delta-1)/2)^2 \int_{G_j}dx\,d_j(x)^{\delta-2} |\varphi(x)|^2
\end{eqnarray*}
for all $\varphi\in C_c^1(G_j)$.
 If, however, one tries to make a similar argument for the boundary layer $\Gamma_{\!\!j,r}$ adjacent to a facial set  $S_j$ then one has to consider the effect of the boundary terms.
There is no problem with the subsets of the boundary of $\Gamma_{\!\!j,r}$ in the interior or exterior faces if $\varphi\in C_c^1(\Gamma_{\!\!r})$ since $\varphi=0$ on these subsets.
But the interfaces $I_{\!j,k}$, $k\neq j$, cause complications  as $\varphi$ cannot be assumed to vanish  on these subsets.

The general situation relies on a multi-dimensional version of (\ref{ehcl1.0}).
If  $\varphi\in C_c^1(\Gamma_{\!\!r})$ then 
\begin{eqnarray*}
0&\leq&\int_{\Gamma_{\!\!j,r}} dx\,\Big|\,d_\Gamma^{\,\delta/2}(\nabla\varphi) +\lambda\,d_\Gamma^{\,\delta/2-1}(\nabla d_\Gamma)\,\varphi\,\Big|^2\\[5pt]
&=&\int_{\Gamma_{\!\!j,r}} dx\,\Big(d_\Gamma^{\,\delta}\,|\nabla\varphi|^2+\lambda^2\,d_\Gamma^{\,\delta-2}\varphi^2+
\lambda\,\delta^{-1}\,(\nabla d_\Gamma^{\,\delta}).(\nabla\varphi^2)\Big)
\\[5pt]
&=&\int_{\Gamma_{\!\!j,r}} dx\,\Big(d_\Gamma^{\,\delta}\,|\nabla\varphi|^2+\lambda^2\,d_\Gamma^{\,\delta-2}\varphi^2-
\lambda\,\delta^{-1}\,(\nabla^2 d_\Gamma^{\,\delta})\,\varphi^2+\lambda\,\delta^{-1}\,\divv(  \varphi^2 \,\nabla d_\Gamma^{\,\delta})\Big)
\end{eqnarray*}
for all $\lambda\in\Ri$ where we have used the  identities $|\nabla d_\Gamma|=1$ on $\Gamma_{\!\!r}$ and
\[
(\nabla d_\Gamma^{\,\delta}).(\nabla\varphi^2)=-(\nabla^2 d_\Gamma^{\,\delta})\,\varphi^2+\divv(  \varphi^2 \,\nabla d_\Gamma^{\,\delta})
\;.
\]
The first of these follows since each $x\in \Gamma_{\!\!j,r} $ has a unique nearest point in $\pi_j$.
The second is the product formula for differentiation. But now one has
\[
(\nabla^2 d_\Gamma^{\,\delta})=\delta(\delta-1)\, d_\Gamma^{\,\delta-2}|\nabla d_\Gamma|^2+\delta\, d_\Gamma^{\,\delta-1}(\nabla^2d_\Gamma)\leq \delta(\delta-1)\, d_\Gamma^{\,\delta-2}
\]
since $-\nabla^2d_\Gamma\geq 0$  by convexity.
These manipulations are all justified  within the integrals because $\varphi\in C_c^1(\Gamma_{\!\!r})$.
Therefore one obtains
\[
\int_{\Gamma_{\!\!j,r}} dx\,\Big(d_\Gamma^{\,\delta}\,|\nabla\varphi|^2+\lambda^2\,d_\Gamma^{\,\delta-2}\varphi^2-
\lambda\,(\delta-1)\,d_\Gamma^{\,\delta-2}\,\varphi^2-\lambda\,\delta^{-1}\,\divv(  \varphi^2 \,\nabla d_\Gamma^{\,\delta})\Big)\geq0
\]
for all $\lambda\geq 0$.
Hence choosing $\lambda=(\delta-1)/2$, which is positive since $\delta>1$, one obtains
\[
\int_{\Gamma_{\!\!j,r}}dx\,d_\Gamma^{\,\delta}\,|\nabla\varphi|^2\geq ((\delta-1)/2)^2\int_{\Gamma_{\!\!j,r}} dx\,d_\Gamma^{\,\delta-2}\,\varphi^2
-((\delta-1)/2\delta)\int_{\Gamma_{\!\!j,r}} dx\,\divv(  \varphi^2 \,\nabla d_\Gamma^{\,\delta})
\]
for all $\varphi\in C_c^1(\Gamma_{\!\!r})$.
Then, however, by the divergence theorem
\[
\delta^{-1}\int_{\Gamma_{\!\!j,r}} dx\,\divv(  \varphi^2 \,\nabla d_\Gamma^{\,\delta})
=\delta^{-1}\int_{\partial\Gamma_{\!\!j,r}} dA\,(n.\nabla d_\Gamma^{\,\delta})\,\varphi^2
=\int_{\partial\Gamma_{\!\!j,r}} dA\,(n.\nabla d_\Gamma)\,(d_\Gamma^{\,\delta-1}\,\varphi^2)
\]
where $dA$ denotes the surface measure and $n$ is the  outward normal, i.e.\ $n.\nabla=\partial/\partial n$ is the normal derivative
at the boundary.
But since $\varphi\in C_c^1(\Gamma_{\!\!r})$ it is zero on the interior and the exterior sections of the boundary of $\Gamma_{\!\!j,r}$ the surface integral
only gives a  possibly non-zero contribution on the interfaces $I_{\!j,k}$ with $k\neq j$.
Therefore combination of these observations gives the Hardy type inequality
\[
\int_{\Gamma_{\!\!j,r}} dx\,d_\Gamma^{\,\delta}\,|\nabla\varphi|^2\geq ((\delta-1)/2)^2\int_{\Gamma_{\!\!j,r}} dx\,d_\Gamma^{\,\delta-2}\,\varphi^2
-((\delta-1)/2)\sum_{k\neq j}\int_{I_{\!j,k}} dA\,(\partial d_\Gamma/\partial n)\,(d_\Gamma^{\,\delta-1}\,\varphi^2)
\]
for all $\varphi\in C_c^1(\Gamma_{\!\!r})$ and $\delta>1$.
This is the multi-dimensional version of the one-dimensional inequality (\ref{ehcl1.0}).

Finally, since $d_\Gamma(x)=d_j(x)=d_k(x)$ for $x\in I_{\!j,k}=I_{\!k,j}$ one has
\[
\int_{I_{\!j,k}} dA\,(\partial d_\Gamma/\partial n)\,(d_\Gamma^{\,\delta-1}\,\varphi^2)+\int_{I_{k,j}} dA\,(\partial d_\Gamma/\partial n)\,(d_\Gamma^{\,\delta-1}\,\varphi^2)
=0
\;.
\]
The flow across the interface $I_{\!j,k}$ from $\Gamma_{\!\!j,r}$ into $\Gamma_{\!\!k,r}$ is cancelled by the flow into $\Gamma_{\!\!j,r}$ from $\Gamma_{\!\!k,r}$ for each pair $j\neq k$.
Therefore 
\begin{eqnarray*}
\int_{\Gamma_{\!\!r}}d_\Gamma^{\,\delta}\,|\nabla\varphi|^2&=&\sum_{j=1}^n\int_{\Gamma_{\!\!j,r}}dx\,d_\Gamma^{\,\delta}\,|\nabla\varphi|^2\\[-5pt]
&\geq&((\delta-1)/2)^2\sum_{j=1}^n\int_{\Gamma_{\!\!j,r}} dx\,d_\Gamma^{\,\delta-2}\,\varphi^2=((\delta-1)/2)^2\int_{\Gamma_{\!\!r}} dx\,d_\Gamma^{\,\delta-2}\,\varphi^2
\end{eqnarray*}
for all $\varphi\in C_c^1(\Gamma_{\!\!r})$.
Thus the boundary Hardy inequality is satisfied  and $a_\delta(\Gamma_{\!\!r})\leq 2/(\delta-1)$ for all $\delta>1$ and all small $r$.
In particular $a_\delta(\Gamma)\leq 2/(\delta-1)$.
\hfill$\Box$

\bigskip

The conclusion of the proposition for convex polytopes is the first step in the proof for general convex domains $\Omega$.
The remainder of the proof involves approximation of $\Omega$ by an increasing family of convex polytopes $T_{\!k}$.
In the sequel $\Gamma$  denotes the boundary of $\Omega$, as in preceding sections,  $\Gamma_{\!\!k}$ 
denotes the boundary of $T_{\!k}$ and $\Gamma_{\!\!k,r}$ denotes the boundary layer of $T_{\!k}$ measured with respect
to the distance $d_{\Gamma_{\!\!k}}$ to the boundary of $T_{\!k}$.
The notation $T$ is used in place of $S$ to avoid confusion with the foregoing discussion in which  $S_k$ denotes a facial set of  the polytope $S$.

\bigskip
\noindent{\bf Proof of Theorem~\ref{tc5}}$\;$
First we  assume that the convex domain $\Omega$ is bounded.
Then, by convexity,  there exists a family of bounded convex polytopes $T_{\!k}\subset \Omega$
which form an increasing  approximating sequence in the following  sense.
First $T_{\!k}\subset T_{\!l}$ for all $l> k$.
Secondly, for each compact subset $K\subset \Omega$ there is a $k$ such that 
$K\subset T_{\!l}\,$ for all $l>k$.
Thirdly $d_{\Gamma_{\!l}}(x)\to d_\Gamma(x)$  as $l\to \infty$ for all $x\in K$.
This is a slight variation of the definition of a normal approximating sequence introduced in \cite{MMP} and \cite{BrM}.

 Secondly fix $\varphi\in C_c^1(\Gamma_{\!\!r})$.
 Let $K=\supp\varphi$.
 Then by the compactness assumption one may choose $k$ such that $\supp\varphi\subset T_{\!l}$
 for all $l>k$.
 But since $T_{\!k}\subset T_{\!l}\subset \Omega$ one has $d_{\Gamma_{\!l}}(x)\leq d_\Gamma(x)$ for all $x\in K$.
Therefore  $\supp\varphi\subset \Gamma_{\!l,r}$ for all $l>k$.
 Consequently,
 \begin{eqnarray*}
 \int_{\Gamma_{\!\!r}} d_\Gamma^{\,\delta}\,|(\nabla\varphi)|^2&\geq&
  \int_{\Gamma_{\!l,r}} d_{\Gamma_{\!l}}^{\,\delta}\,|(\nabla\varphi)|^2\\
  &\geq&((\delta-1)/2)^2 \int_{\Gamma_{\!l,r}} d_{\Gamma_{\!l}}^{\,\delta-2}\,|\varphi|^2\geq ((\delta-1)/2)^2 \int_{\Gamma_{\!\!r}} d_{\Gamma_{\!l}}^{\,\delta-2}\,|\varphi|^2
  \end{eqnarray*}
  for all $l>k$ where, in the second step,  we have applied Proposition~\ref{pcnotes1} to the polytopes $T_{\!l}$.
  Now if $\delta\in\langle1,2\,]$ then $d_{\Gamma_{\!l}}^{\,\delta-2}\geq d_{\Gamma}^{\,\delta-2}$, because $d_{\Gamma_{\!l}}\leq d_{\Gamma}$ on $\supp\varphi$,
  and one obtains the weighted Hardy inequality on $\Gamma_{\!\!r}$.
  Alternatively, if $\delta>2$ then the right hand integral is bounded uniformly in $l$.
   In addition $d_\Gamma^{\,\delta-2}|\varphi|^2$ is bounded on $\Gamma_{\!\!r}$.
 Then it follows from the Lebesgue dominated convergence theorem and  the condition $d_{\Gamma_{\!l}}\to d_{\Gamma}$ on $\supp\varphi$ that 
\begin{equation}
  \int_{\Gamma_{\!\!r}} d_\Gamma^{\,\delta}\,|(\nabla\varphi)|^2\geq ((\delta-1)/2)^2 \int_{\Gamma_{\!\!r}} d_{\Gamma}^{\,\delta-2}\,|\varphi|^2
  \label{ehc1.01}
  \end{equation}
  for one and hence all $\varphi\in C_c^1(\Gamma_{\!\!r})$.
  Therefore $a_\delta(\Gamma_{\!\!r})\leq 2/(\delta-1)$ for all $\delta>1$ and all small $r>0$.
  Hence $a_\delta(\Gamma)\leq 2/(\delta-1)$.
  This completes the proof for bounded $\Omega$.
  
  Finally if $\Omega$ is unbounded fix a point $z\in\Gamma$ and set $\widetilde \Omega=\Omega\cap B(z\,;R)$.
  Let $\widetilde \Gamma$ denote the boundary of  $\widetilde \Omega$.
  Then if $\varphi\in C_c^1(\Gamma_{\!\!r})$ one may choose $z, R$ such that $\supp\varphi\subset {\widetilde \Gamma}_{\!\!r}$.
  Hence  $\varphi\in C_c^1({\widetilde \Gamma}_{\!\!r})$ and $d_{\tilde\Gamma}(x)=d_\Gamma(x)$ for all $x\in \supp\varphi$.
  Consequently (\ref{ehc1.01}) follows directly from the result for bounded domains.
  It suffices to interpret the inequality on the bounded domain $\widetilde \Omega$.
  \hfill$\Box$

\section{The complement of convex sets}\label{S6}

In this section we examine the Hardy inequality for domains which are the complement of convex sets.
Specifically $\Omega=\Ri^d\backslash K$ where $K$ is a closed convex subset of $\Ri^d$.
We assume throughout that $K$ is non-trivial and  $d$-dimensional.
Then the Hausdorff dimension of the boundary $\Gamma$ of $\Omega$ is $d-1$.
Again the dichotomy of Theorem~\ref{thc3.1}  corresponds to the two cases $\delta\in[0,1\rangle$ and
$\delta>1$.
In contrast to the case of convex domains the $\delta>1$ regime is well understood but little attention
has been paid to the  $\delta\in[0,1\rangle$ case which in fact deviates from the behaviour encountered in the earlier sections.
First we summarize the relevant properties for $\delta>1$ in the following proposition which follows directly from Theorem~1.1 in \cite{Rob13}
although  we give a short independent proof based on  the arguments of Section~\ref{S5}.

\begin{prop}\label{phc6.10}
If $\Omega=\Ri^d\backslash K$ is the complement of a non-trivial, closed, $d$-dimensional convex subset $K$ of $\Ri^d$ and $\delta>1$
then
the weighted  Hardy inequality $(\ref{ehc1.1})$ is satisfied for all $\varphi\in C_c^1(\Omega)$ and  $a_\delta(\Omega)=a_\delta(\Gamma)=2/(\delta-1)$.
\end{prop}
\proof\
The proof in  \cite{Rob13} depends on the observation that each point $x\in \Omega$ has a unique nearest point $n(x)\in \Gamma$, by Motzkin's theorem.
This implies that $d_\Gamma$ is differentiable and $|(\nabla d_\Gamma)(x)|=1$ for all $x\in \Omega$ (see, for example, \cite{BEL}, Chapter~2).
Moreover, $d_\Gamma$ is convex on convex  subsets of $\Omega$ and so $\nabla^2d_\Gamma$ is a positive measure.
But then one has the following shorter proof based on the foregoing discussion.

First since $\delta>1$ and the integral of  $(\nabla d_\Gamma).(\nabla d_\Gamma^{\,\delta-1}\varphi^2)$ over $\Omega$ is negative since  $\nabla^2d_\Gamma$ is positive
it follows from the  identity  (\ref{ihc1}) that 
\[
0\leq (\delta-1)\,\|d_\Gamma^{\,\delta/2-1}\varphi\|_2^2\leq \Big|\int_{\Omega}\,d_\Gamma^{\,\delta-1}((\nabla d_\Gamma).(\nabla\varphi^2))\Big|
\leq 2\,\|d_\Gamma^{\,\delta}(\nabla\varphi)\|_2\,\|d_\Gamma^{\,\delta/2-1}\varphi\|_2
\]
for all   $\varphi\in C_c^1(\Omega)$. 
Hence the Hardy inequality is valid on $\Omega$ for all $\delta>1$ and one  has the upper bound $a_\delta(\Omega)\leq 2/(\delta-1)$.
Moreover,  the lower bound $2/(\delta-1)\leq a_\delta(\Gamma)\leq a_\delta(\Omega)$ 
follows from Proposition~\ref{prhc3.1}. 
\hfill$\Box$

\bigskip
The situation with  $\delta\in[0,1\rangle$ is more complicated.
But the Hardy constant $a_\delta(\Gamma)=2/(1-\delta)$ for all $\delta\in[0,1\rangle$  if and only if $a_0(\Gamma)\leq 2$
by the discussion at the beginning of Section~\ref{S4} and the lower bound of Proposition~\ref{prhc3.1}.
Since $a_0(\Gamma)=b_0(\Gamma)=b_0(\Omega)$, by Proposition~\ref{phc2.1},
this latter condition is also  equivalent to   $b_0(\Omega)\leq 2$.
Nevertheless, if $d\geq2$  then there are $K$ such that $b_0(\Omega)>2$.
In fact $K$ can be a bounded polytope.

\begin{prop}\label{phc6.2}
Let $S$ be a bounded, closed, convex, $d$-dimensional polytope and set $\Omega=\Ri^d\backslash S$.
Further let $\{\pi_j\}_{\{1\leq j\leq n\}}$ denote the hyperplanes containing the faces of $S$ and
$\{\alpha_{jk}\}_{\{1\leq j\leq n\}}$ the dihedral angles corresponding to the pairs $\pi_j, \pi_k$ with $j\neq k$.

Then there is an $\alpha_c\in\langle0,\pi\rangle$ such that $a_0(\Gamma)=b_0(\Gamma)=b_0(\Omega)>2$
 if $\alpha_{jk}<\alpha_c$ for some choice of $j,k$.
\end{prop}
\proof\
Again let $\{\Gamma_{\!\!j}\}_{\{1\leq j\leq n\}}$ denote the open faces of $S$.
Thus $\Gamma=\bigcup^{\,n}_{j=1}{\overline\Gamma}_{\!\!j}$.
Further let $G_j$ denote the half space  such that $(\partial G_j\cap S)={\overline\Gamma}_{\!\!j}$.
In contrast to the situation in Section~\ref{S5} the polytope is  now in the complement  of each of the half-spaces.
Moreover, $\Omega=\bigcup^{\,n}_{j=1}G_j$.
Now to estimate $b_0(\Omega)$ we use the identification $b_0(\Omega)=b_0(\Gamma)=a_0(\Gamma)$ of Proposition~\ref{phc2.1}
together with the identification of $a_0(x)$ and $b_0(x)$ for $x\in\Gamma$ given by Theorem~\ref{thc2.1}.
The latter identifies $a_0(\Gamma)$ as the supremum of the local function $x\in \Gamma\mapsto a_0(x)=\inf\{a_0(\Gamma\cap U): U\subset \Ri^d, \,U\ni x\}$.
Since $b_0(\Gamma\cap U)=a_0(\Gamma\cap U)$ by Proposition~\ref{phc2.1} one can equally well identify $a_0(x)$ with the weak local constant $b_0(x)$.
Hence to deduce that $b_0(\Gamma)>2$ it suffices to identify one point $x\in \Gamma$ such that $b_0(x)>2$.
We choose a point in the interior of the $(d-2)$-dimensional  edge  formed by the intersection ${\overline\Gamma}_{\!\!j}\cap{\overline\Gamma}_{\!\!k}$ of the closures of two of the faces
and set $\alpha_{j,k}=\alpha$ for simplicity.

The point  $x$ is not only a boundary point of $\Omega$ but it is also a boundary point of the domain $\Omega_{jk}=G_j\cup\, G_k$ formed by the half-spaces containing the faces. 
Therefore, since  $b_0(x)$ is purely local, it can be computed with respect to the domain $\Omega_{jk}$.
Now  $\Omega_{jk}=\Ri^d\backslash S_{jk}$ where $S_{jk}$ is a closed, unbounded, convex polytope  with two faces contained in the faces of $G_j$ and $G_k$, respectively.
One can then  choose coordinates such that $\Ri^d=\Ri^2\times \Ri^{d-2}$ and $S_{jk }=T_{\!jk}\times \Ri^{d-2}$
where $T_{\!jk}\subset \Ri^2$ is a $2$-dimensional unbounded polytope  given by 
\[
T_{\!jk}=\{ y=re^{i\theta}: |y|=r\geq0, \mbox{ and } |\theta|\leq \alpha/2<\pi/2\}\;,
\]
i.e., it  is a pointed cone with apex at the origin and two faces $re^{\pm i\alpha/2}$.
The angle $\alpha$ is the dihedral angle between the two half-spaces  $G_j$ and $G_k$
and the condition $\alpha<\pi$ follows from the convexity of $S$.
Now one can reduce the problem of estimating $b_0(\Omega_{jk})$ to a $2$-dimensional problem.

\begin{lemma}\label{lhc6.1}
Let ${\widetilde \Omega}_{jk}=\Ri^2\backslash T_{jk}$.
Then the  weak Hardy constant  $b_0(\Omega_{jk })$ is larger or equal to the weak Hardy constant $ b_0({\widetilde \Omega}_{jk})$.
\end{lemma}
\proof\
It follows by construction that $\Omega_{jk}={\widetilde \Omega}_{jk}\times \Ri^{d-2}$.
So if $x\in\Omega_{jk}$ then $x=(y,z)$ with $y\in {\widetilde \Omega}_{jk}$ and $z\in \Ri^{d-2}$.
Moreover, the distance $d_\Gamma(x)$ to the boundary of $\Omega_{jk}$ is equal to the distance $d_\Gamma(y)$ of $y$ to the boundary of ${\widetilde \Omega}_{jk}$.
It is independent of $z$.
Now the weak Hardy constant  $b_0(\Omega_{jk})$ is the infimum of the $b>0$ for which there is a $c\geq 0$ such that
\[
\|d_\Gamma^{\,-1}\varphi\|_2\leq b\,\|\nabla\!\varphi\|_2+c\,\|\varphi\|_2
\]
for all  $\varphi\in C_c^1(\Omega_{jk})$.
But the inequality is valid for $\varphi$ which are products of $C_c^1$-functions in the two variables.
So setting $\varphi=\psi\chi$, with $\psi\in C_c^1({\widetilde\Omega}_{jk})$ and $\chi\in C_c^1(\Ri^{d-2})$, and then dividing by $\|\chi\|_2$ one obtains
\begin{eqnarray*}
\|d_\Gamma^{\,-1}\psi\|_2&\leq& b\,\|\nabla\!\psi\|_2+{\tilde c}\,\|\psi\|_2
\end{eqnarray*}
for all $\psi\in C_c^1({\widetilde\Omega}_{jk})$ where ${\tilde c}=c+b\,\|\nabla\! \chi\|_2/\|\chi\|_2$.
Therefore  $ b_0({\widetilde \Omega}_{jk})\leq b$.
Since $b_0({\Omega}_{jk})$ is the infimum over the possible choices of $b$  it follows that 
$ b_0({\widetilde \Omega}_{jk})\leq   b_0({\Omega}_{jk})$.\hfill$\Box$

\bigskip

Now if $x\in {\overline\Gamma}_{\!\!j}\cap{\overline\Gamma}_{\!\!k}=\pi_j\cap\pi_k$ then $x=(0,z)$ for some $z\in\Ri^{d-2}$.
Therefore Lemma~\ref{lhc6.1}  establishes that $b_0(x)$  is larger  or  equal to the weak Hardy constant corresponding to ${\widetilde\Omega}_{jk}$
at the apex of the cone $T_{\!jk}$, i.e.\ at the origin in $\Ri^2$.
Adopting  the notation $b_{0,\alpha}$ for this  $2$-dimensional Hardy constant one concludes that $b_0(x)\geq b_{0,\alpha}(0)$.
But it follows from Davies' analysis of two-dimensional  sectors,  \cite{Dav17} Section~4, that  there is a critical $\alpha_c\in\langle0,\pi\rangle$ such that $b_{0,\alpha}(0)>2$ if and only if $\alpha<\alpha_c$. 
Therefore $b_0(x)>2$ if   $\alpha<\alpha_c$, i.e.\
 $b_0(\Omega)>2$ if any one of the dihedral angles defined by  pairs of  faces $\Gamma_{\!\!j}, \Gamma_{\!\!k}$ of $S$  is less than $\alpha_c$.
\hfill$\Box$

\bigskip

We have adopted two different notational conventions to Davies so we comment on the details of the application of his results.
First his weak Hardy constant is the square of our constant.
Secondly, Davies resolved the $2$-dimensional problem of conic sectors by rephrasing the problem as a Hardy inequality in terms of the angular coordinate $\theta$.
In particular he considered the exterior angle $\beta$ of the cone instead of the interior angle $\alpha$ that we have used.
We preferred the latter since it  corresponds to the dihedral angle between the half-spaces.
For comparison $\beta\in[\pi, 2\pi\rangle$ and $\beta-\pi=\pi-\alpha$.
Davies  then  deduced from the angular Hardy inequality, by numerical analysis,  that the weak Hardy constant assumes the standard value if and only if $\beta\leq \beta_c$
where  $\beta_c$ is approximately $4.856$ radians, or $1.546\, \pi$.
Thus with our conventions the corresponding critical angle $\alpha_c$ is approximately $0.454\,\pi$ or $81.77^\circ$ and $b_0(x)=2$ for all $x\in \Gamma$ if and only if $\alpha\geq \alpha_c$.
Subsequently, Tidblom \cite{Tid}  solved the angular Hardy inequality in terms of Legendre functions of the first kind and found that
\[
\beta_c-\pi=4\,\arctan(2\Gamma(3/4)/\Gamma(1/4))^2=\pi-\alpha_c
\]
although Tidblom also used a different convention.
(An alternative derivation of this expression and a formula for the  the Hardy constant for $\beta >\beta_c$ can be found in  \cite{BT1} \cite{BT2}.)
The significance of this expression is not clear but it is certainly the first step in understanding the geometric factors that govern the 
Hardy inequality in higher-dimensions.

The strength of the foregoing argument is that it draws a multi-dimension conclusion from a two-dimensional result.
But that is also its weakness. 
The reduction lemma, Lemma~\ref{lhc6.1}, ensures that any ensuing condition for anomalous values of the Hardy constants will only be a necessary condition.
The calculation aims to gain an upper bound on the local weak Hardy constant $b_0(x)$ but the bound is restricted to points $x$ in the  $(d-1)$-dimensional faces of the polytope $S$
or in the $(d-2)$-dimensional edges.
But the value of $b_0$ is expected to increase as one passes to lower dimensional edges and to be maximal at the vertices of the polytope.
Although this intuition does not have any substantial, quantitive, foundation.
Nevertheless one could aspire to tackling the problem in a systematic manner by passing successively to the examination of points  in
$(d-3)$-dimensional edges contained in the intersection of three of the hyperplanes $\pi_j$ and then to the  $(d-4)$-dimensional edges et cetera.
Each of these problems is reduced to a lower-order problem  just as the proof of Proposition~\ref{phc6.2} was reduced to a two-dimensional problem.
It should be feasible to resolve the analogous three-dimensional problem by numerical analysis or at least to derive some further information concerning the Hardy constant.

Two illustrations of Proposition~\ref{phc6.10} are given by the following examples.

\begin{exam}\label{ex6.1} Let $\Omega=\Ri^d\backslash C_d$ where $C_d$ is a closed, $d$-dimensional, rectangle.
Then all dihedral angles are equal to $\alpha=\pi/2>\alpha_c$.
Hence $a_0(\Gamma)=b_0(\Gamma)=b_0(\Omega)=2$ for all dimensions. 
Then $a_\delta(\Gamma)=b_\delta(\Gamma)=b_\delta(\Omega)=2/(1-\delta)$ for all $\delta\in[0,1\rangle$ by Proposition~\ref{phc2.1}
and the argument at the beginning of Section~\ref{S5}.
In combination with Proposition~\ref{phc6.10} one deduces that the standard situation is valid for all $\delta\geq0$ with $\delta\neq1$.
A similar conclusion is valid if $C_d$ is a $d$-dimensional rhomboid as long as the dihedral angles remain larger than~$\alpha_c$.
\end{exam}

\begin{exam}\label{ex6.2}
Let  $\Omega=\Ri^d\backslash S_d$ where $S_d$ is the regular $d$-simplex.
Then there is a unique dihedral angle $\alpha= \arccos(1/d)$  between the various faces (see, for example, \cite{PaW}).
Then $a_0(\Gamma)=b_0(\Gamma)=b_0(\Omega)>2$ if $d\geq 7$.
In fact one has the approximate value $\cos\alpha_c=0.1431$ and this is slightly larger than $1/7$ but smaller than $1/6$.
Moreover, $a_0(\Gamma)=b_0(\Gamma)=b_0(\Omega)=2$ if $d=1$ or $2$.
The one-dimensional case is clear and the two-dimensional case follows from \cite{Dav17}, Theorem~4.1.
The  situation with $d=3,4,5, 6$ is unclear.
In each of these case $\alpha>\alpha_c$ but this  criterion for a $(d-2)$-dimensional edge is insufficient to draw a conclusion on the value of the Hardy
constant.
If $x_j$ is an arbitrary point in a $(d-j)$-dimensional edge then $b_0(\Omega)=\max_{\{1\leq j\leq d\}}(b_0(x_j))$ 
and it is expected that the maximum is $b_0(x_d)$, i.e. the maximum is attained at a vertex of the simplex.
It would be interesting to verify this, numerically at least,  in the simplest case $d=3$.
\end{exam}

 \section{Operator implications}\label{S7}

The prime motivation for the the preceding analysis stemmed from earlier work  \cite{Rob16} on the self-adjointness of symmetric degenerate elliptic operators  of the form $H=-\divv(C\,\nabla)$
on $L_2(\Omega)$.
We briefly summarize the essential definitions of the previous paper and then derive self-adjointness criteria on $C^{1,1}$-domains, convex domains etc.\  based on the foregoing analysis.

Let $C=(c_{kl})$ be a  strictly positive, symmetric, $d\times d$-matrix with  $c_{kl}$ real, Lipschitz continuous, functions 
which  resembles the diagonal matrix $c\,d_\Gamma^{\,\delta}I$.
Specifically,  we assume 
 \begin{equation}
\textstyle{\inf_{r\in\langle0,r_0]}}\;\textstyle{\sup_{x\in\Gamma_{\!\!r}}}\|(C\,d_\Gamma^{\,-\delta})(x)-c(x)I\|=0
\label{esa5.1}
 \end{equation}
 for some $r_0>0$
 where $c$ is  a  bounded Lipschitz function satisfying  $\inf_{x\in\Gamma_{\!\!r}}c(x)\geq \mu>0$ and $\delta\geq0$.
Condition~(\ref{esa5.1})  can be interpreted in an obvious way as
\[
\textstyle{\limsup_{d_\Gamma\to0}} \;C(c\,d_\Gamma^{\,\delta}I)^{-1}=I
\;.
\]
Thus in the language of asymptotic analysis  $C$ converges to  $c\,d_\Gamma^{\,\delta} I$ as $d_\Gamma\to0$ (see \cite{DeB}).
The  parameter $\delta$ determines the order of degeneracy at the boundary and  $c$ describes  the boundary  profile of $C$. 
With these assumptions $H$ is defined as a positive symmetric operator on $C_c^\infty(\Omega)$ and we  also use $H$ to
denote its symmetric closure.

The operator $H$ is not necessarily self-adjoint and an obvious  problem is to obtain necessary and sufficient conditions which
ensure self-adjointness.
A partial step in this direction is given by Theorem~1.1 in \cite{LR} which deals with the problem of Markov uniqueness, i.e.\
the existence of a unique self-adjoint extension of $H$ which generates a Markov semigroup on $L_2(\Omega)$.
Under the assumption that $\Omega$ is a uniform domain with an Ahlfors regular boundary it was established that $H$
is Markov unique if and only if $\delta\geq 2-(d-d_{\!H})$. 
Thus this condition is necessary for self-adjointness but not sufficient.
On the other hand if $\delta\geq 2$ then $H$ is self-adjoint (with no uniformity or regularity restrictions  on the domain). 
This is established in Corollary~2.5 of \cite{Rob16} although the result has many precedents.
Therefore the condition $\delta\geq 2$ is sufficient for self-adjointness and it remains to consider the range of $\delta\in [2-(d-d_{\!H}), 2\rangle$.

\begin{thm}\label{thc7.1}
Assume $\Omega$ is a uniform domain with an Ahlfors regular boundary and that $\Omega$ is either a $C^{1,1}$-domain, or a convex domain or  the complement
of a closed convex subset of $\Ri^d$. 
Further assume the coefficients of $H$ satisfy the boundary condition $(\ref{esa5.1})$.
Then  the condition $\delta>3/2$ is sufficient for self-adjointness of $H$.

If, in addition,  
\begin{equation}
\textstyle{\sup_{x\in\Gamma_{\!\!r_0}}}|(\divv (Cd^{\,-\delta})).(\nabla d_\Gamma)(x)|<\infty
\;. \label{eun5.21}
 \end{equation}
Then the condition $\delta\geq 3/2$ is necessary for self-adjointness of $H$.
\end{thm}
\proof\
The first statement  is a direct corollary of Theorem~5.2 in \cite{Rob16} together with the results of the preceding sections.
Theorem~5.2 states, in the current notation, that if $\delta\in [1,2\rangle$ then $a_\delta(\Gamma)>2/(2-\delta)$ is sufficient for self-adjointness of $H$.
But in all three cases of the theorem one now has $a_\delta(\Gamma)=2/(\delta-1)$ and the sufficient condition reduces to $\delta>3/2$.

The second statement of the theorem is a generalisation of a
 result  derived for $C^2$-domains in \cite{Rob15}, Theorem3.2.
The proof 
 for the current situations is identical. 
 It depends on demonstrating that if $\delta<3/2$ then there is a non-zero $\varphi$ in the domain of the adjoint $H^*$ of $H$
 but $\varphi$ is not in the quadratic form domain. This implies that $H$ is not self-adjoint.
\hfill$\Box$

\bigskip

Theorem~\ref{thc7.1} also extends to domains which  are the complement of  a countable family
of convex or $C^{1,1}$-domains which are positively separated since the  boundary Hardy inequality  separates into a family of inequalities on the boundaries of 
of each domain  with the same Hardy constant.
This follows since all the crucial estimates involved in the proof are restricted to boundary layers (for further details see \cite{Rob15}).

\newpage

Finally we demonstrate that the boundary Hardy inequality is equivalent to a spectral property of the self-adjoint Friedrichs' extension  of $H$
at least for  $\delta\in[0,2\rangle$.

\begin{thm}\label{thc7.2}
Let $\Omega$ be a general domain and assume the coefficients of  $H=-\divv(C\,\nabla)$
on $L_2(\Omega)$ satisfy the boundary condition $(\ref{esa5.1})$.
Further let  $H_\delta$ denote   the self-adjoint  Friedrichs' extension  of $H$.

 If $\delta\in[0,2\rangle$ then the following conditions are equivalent.
\begin{tabel}
\item\label{phc7.2-1}
The boundary  Hardy inequality $(\ref{ehc1.1})$ is valid on $\Gamma_{\!\!r}$ for some small $r>0$,
\item\label{phc7.2-2}
there is a $\beta>0$ such that $H_\delta-\beta \,c\,d_\Gamma^{\,\delta-2}I$ is lower semibounded on $L_2(\Omega)$.
\end{tabel}
Moreover, if these conditions are satisfied then the supremum over the possible $\beta$ in the second condition  is equal to $a_\delta(\Gamma)^{-2}$.
\end{thm}
\proof\
As a preliminary note that only the values of $c$ on $\Gamma_{\!\!r_0}$  are relevant in the boundary condition (\ref{esa5.1}).
Nevertheless one may assume that $c$ extends to $\Omega$.
Moreover (\ref{esa5.1}) implies that for each $r\in\langle0,r_0]$ there is $\sigma_r\geq 0$ such that $0\leq (c\,d_\Gamma^{\,\delta})I\leq \sigma_rC$ on $\Gamma_{\!\!r}$ and $\sigma_r\to1$ as $r\to0$.

First assume the boundary Hardy inequality is valid on $\Gamma_{\!\!r}$.
Then
\begin{eqnarray*}
\|c^{1/2}d_\Gamma^{\,\delta/2-1}\varphi\|_2&=&\|d_\Gamma^{\,\delta/2-1}(c^{1/2}\varphi)\|_2\\[5pt]
&\leq& a\,\|d_\Gamma^{\,\delta/2}(\nabla (c^{1/2}\varphi))\|_2\leq a\,\|c^{1/2}d_\Gamma^{\,\delta/2}(\nabla \!\varphi)\|_2+a\,\|\nabla c^{1/2}\|_\infty\|d_\Gamma^{\,\delta/2}\varphi\|_2
\end{eqnarray*}
for all $\varphi\in C_c^1(\Gamma_{\!\!r})$.
But  $\|d_\Gamma^{\,\delta/2}\varphi\|_2= r^{\delta/2}\|(d_\Gamma/r)^{\delta/2}\varphi\|_2\leq r^{\delta/2} \|\varphi\|_2$.
Combining these estimates one obtains
\[
\|c^{1/2}d_\Gamma^{\,\delta/2-1}\varphi\|_2\leq a\,\|c^{1/2}d_\Gamma^{\,\delta/2}(\nabla \!\varphi)\|_2+a_r\,\|\varphi\|_2
\]
where $a_r=a\,\|\nabla c^{1/2}\|_\infty\,r^{\delta/2}$. 
Then if $\varepsilon\in\langle 0,1\rangle$ and one chooses $r$ small enough that $a_r<\varepsilon$ one finds
\begin{equation}
\|c^{1/2}d_\Gamma^{\,\delta/2-1}\varphi\|_2\leq (1-\varepsilon)^{-1}a\,\|c^{1/2}d_\Gamma^{\,\delta/2}(\nabla \!\varphi)\|_2\leq (1-\varepsilon)^{-1}a \,\sigma_{r_0}^{1/2}\,h(\varphi)^{1/2}
\label{ehc7.3}
\end{equation}
for all $\varphi\in C_c^1(\Gamma_{\!\!r})$  and $r\in\langle 0,r_0\rangle$ where $h$ is the quadratic form corresponding to $H$.
Since $\varepsilon$ can be arbitrarily small and $\sigma_{r_0}$ arbitrarily close to one the value of the constant in this  modified boundary inequality is essentially  equal to $a$.

Next one obtains a direct analogue of  Proposition~\ref{phc2.1} for the modified inequalities with no essential change in the proof.
Then it follows that the modified boundary inequality is equivalent to a modified version of the weak Hardy inequality on $L_2(\Omega)$ similar to (\ref{ehc2.3}).
Explicitly, there is a $c_{\varepsilon,r_0}>0$ such that
\[
\|c^{1/2}d_\Gamma^{\,\delta/2-1}\varphi\|_2\leq (1-\varepsilon)^{-1}a \,\sigma_{r_0}^{1/2}\,h(\varphi)^{1/2}
\,\|d_\Gamma^{\,\delta/2}(\nabla\!\varphi)\|_2+c_{\varepsilon,r_0}\,\|\varphi\|_2
\]
for all $\varphi\in C_c^1(\Omega)$.
Then squaring this inequality and applying the Cauchy-Schwarz inequality to the right hand side one deduces that for all $\varepsilon_1>0$ there is a $c_{\varepsilon, \varepsilon_1, r_0}>0$ such that
\[
(\varphi, c\,d_\Gamma^{\,\delta-2}\varphi)\leq (1+\varepsilon_1)(1-\varepsilon)^{-2}a^2\sigma_{r_0}\,h(\varphi)+c_{\varepsilon, \varepsilon_1, r_0}\,(\varphi,\varphi)
\]
for all $\varphi\in C_c^2(\Omega)$.
This inequality then extends to the Friedrichs' extension $H_\delta$ by closure.
Therefore one has a family of operator  bounds
\[
\sigma_{r_0}\,H_\delta\geq \beta_{\varepsilon, \varepsilon_1}\,c\,d_\Gamma^{\,\delta-2}\,I-\gamma_{\varepsilon, \varepsilon_1, r_0} \,I
\]
where $\beta_{\varepsilon, \varepsilon_1}\leq a^{-2}$.
Moreover, $\beta_{\varepsilon, \varepsilon_1}\to a^{-2}$ as $\varepsilon_1, \varepsilon\to0$ 
and  $\sigma_{r_0}\to 1$ as $r_0\to0$. One concludes that for each $\beta< a^{-2}$ there is a $\gamma>0$ such that $H_\delta\geq\beta\,c\,d_\Gamma^{\,\delta-2}\,I-\gamma I$.
Thus  one  verifies simultaneously Condition~\ref{phc7.2-2} and the last statement of the theorem.

Since  the implication \ref{phc7.2-2}$\Rightarrow$\ref{phc7.2-1} is straightforward the proof is complete.\hfill$\Box$

\bigskip

The semiboundedness property for $H$ given by Theorem~\ref{thc7.2} is automatically satisfied if $\Omega$ is a uniform domain with Ahlfors regular boundary
 if $\delta\geq0$ with the exception of $\delta=2-(d-d_{\!H})$ by Theorem~\ref{thc3.1}.
 Hence the equivalent conditions of Theorem~\ref{thc7.2} are satisfied for all $\delta\in[0,2\rangle$ with the exception of the one special value.
 A further improvement occurs in the case of convex or $C^{1,1}$-domains.
 Then $d_{\!H}=d-1$ and so the boundary Hardy inequality is satisfied for all $\delta\in[0,2\rangle$ with the exception of $\delta=1$
 and $a_\delta=2/|1-\delta|$. 
 Thus $H_\delta-\beta\,c\,d_\Gamma^{\,\delta-2}I$ is lower semibounded for all $\beta<(1-\delta)^2/4$.

 \section*{Acknowledgements}
Again I am indebted to Juha Lehrb{\"a}ck  for a wealth of information concerning the Hardy inequality. 
In particular for drawing my attention to the paper by  Haj{\l}asz involving the boundary inequality and correcting various infelicities in 
my earlier draft of Section~\ref{S3}.

\end{document}